\documentclass[a4paper,11pt,reqno,noindent]{amsart}
\usepackage[centertags]{amsmath}
\usepackage{amsfonts,amssymb,amsthm,dsfont,cases,amscd,esint,enumerate}
\usepackage[T1]{fontenc}
\usepackage[english]{babel}
\usepackage[applemac]{inputenc}
\usepackage{newlfont}
\usepackage{color}
\usepackage[body={15cm,21.5cm},centering]{geometry} 
\usepackage{fancyhdr}
\usepackage{enumerate}

\theoremstyle{plain}
\newtheorem{theor0}{Theorem}
\newenvironment{theor}
  {\pushQED{\qed}\begin{theor0}}
  {\popQED\end{theor0}}
\newtheorem{prop10}{Proposition}
\newenvironment{prop1}
  {\pushQED{\qed}\begin{prop10}}
  {\popQED\end{prop10}}
  \newtheorem{cor0}{Corollary}
\newenvironment{cor}
  {\pushQED{\qed}\begin{cor0}}
  {\popQED\end{cor0}}
\newtheorem{lem0}{Lemma}[section]
\newenvironment{lem}
  {\pushQED{\qed}\begin{lem0}}
  {\popQED\end{lem0}}
\newtheorem{prop0}[lem0]{Proposition}
\newenvironment{prop}
  {\pushQED{\qed}\begin{prop0}}
  {\popQED\end{prop0}}
\theoremstyle{definition}
\newtheorem{defin0}[lem0]{Definition}

\newtheorem{rems0}[lem0]{Remarks}

\newtheorem{rem0}[lem0]{Remark}
\newenvironment{rem}
  {\pushQED{\qed}\begin{rem0}}
  {\popQED\end{rem0}}

\numberwithin{equation}{section}

\newcommand{\e}{\varepsilon}

\newcommand{\calP}{\mathcal{P}}

\newcommand{\R}{\mathbb R}

\newcommand{\G}{G}
\newcommand{\F}{F}

\newcommand{\Id}{\operatorname{Id}}

\newcommand{\ee}{e}
\newcommand{\Aa}{\boldsymbol a}

\newcommand{\Ld}{\operatorname{L}}

\newcommand{\step}[1]{\noindent \textit{Step} #1.}
\newcommand{\substep}[1]{\noindent \textit{Substep} #1.}
\newcommand{\Pm}{\mathbb{P}}

\newcommand{\expec}[1]{\mathbb{E}\left[ #1 \right]}
\newcommand{\expecm}[1]{\mathbb{E}\big[ #1 \big]}
\newcommand{\var}[1]{\mathrm{Var}\left[#1\right]}

\newcommand{\ent}[1]{\mathrm{Ent}\!\left[#1\right]}

\newcommand{\calC}{\mathcal{C}}

\def\Aa{\boldsymbol a}
\newcommand{\bb}{\boldsymbol b}

\newcommand{\Norm}[1]{|\!|\!|#1|\!|\!|}

\newcommand{\fun}{\mathrm{fct}}

\usepackage[colorlinks,citecolor=black,linkcolor=black,urlcolor=black]{hyperref}

\title[Robustness of the pathwise structure of fluctuations]{Robustness of the pathwise structure of fluctuations in stochastic homogenization}
\author[M. Duerinckx]{Mitia Duerinckx}
\author[A. Gloria]{Antoine Gloria}
\author[F. Otto]{Felix Otto}
\address[Mitia Duerinckx]{Laboratoire de Mathématique d'Orsay, UMR 8628, Université Paris-Sud, F-91405 Orsay, France \& Universit\'e Libre de Bruxelles, Département de Mathématique, Brussels, Belgium}
\email{mduerinc@ulb.ac.be}
\address[Antoine Gloria]{Sorbonne Universit\'e, CNRS, Universit\'e de Paris, Laboratoire Jacques-Louis Lions (LJLL), F-75005 Paris, France \& Universit\'e Libre de Bruxelles, Département de Mathématique, Brussels, Belgium}\email{antoine.gloria@upmc.fr}
\address[Felix Otto]{Max Planck Institute for Mathematics in the Sciences, Leipzig, Germany}
\email{otto@mis.mpg.de}
\begin{document}
\maketitle

\begin{abstract}
We consider a linear elliptic system in divergence form with random coefficients and study the random fluctuations of large-scale averages of the field and the flux of the solution operator. In the context of the random conductance model, we
developed in a previous work a theory of fluctuations based on the notion of homogenization commutator:
we proved that the two-scale expansion of this special quantity is accurate at leading order in the fluctuation scaling when averaged on large scales (as opposed to the two-scale expansion of the solution operator taken separately) and that the large-scale fluctuations of the field and the flux of the solution operator can be recovered from those of the commutator. This implies that the large-scale fluctuations of the commutator of the corrector drive all other large-scale fluctuations to leading order, which we refer to as the \emph{pathwise structure} of fluctuations in stochastic homogenization.
In the present contribution we extend this result in two directions: we treat \emph{continuum} elliptic (possibly non-symmetric) systems and allow for strongly \emph{correlated} coefficient fields (Gaussian-like with a covariance function that can display an arbitrarily slow algebraic decay at infinity). Our main result shows in this general setting that the two-scale expansion of the homogenization commutator is still accurate to leading order when averaged on large scales, which illustrates the robustness of the pathwise structure of fluctuations.
\end{abstract}

\tableofcontents

\section{Introduction}

This article constitutes the second part of a series that develops a theory of fluctuations in stochastic homogenization of linear elliptic (non-necessarily symmetric) systems in divergence form.
In the first part~\cite{DGO1}, we provided a complete picture of our theory (with optimal error estimates and convergence rates) in the simplified framework of the random conductance model. We proved three main results: the pathwise structure of fluctuations, their asymptotic normality, and the identification of the limiting covariance structure.
In the present contribution, we focus on the fundamental pathwise aspect of the theory, that is, the accuracy of the two-scale expansion for large-scale fluctuations of the so-called homogenization commutator,
and we extend its validity to continuum (non-symmetric) systems with strongly correlated coefficient fields.
More precisely, we cover the general setting of coefficient fields satisfying multiscale functional inequalities as introduced in~\cite{DG1,DG2}, and therefore treat all the models considered in the reference textbook~\cite{Torquato-02} on heterogenous materials. We take this as a sign of the robustness of the pathwise structure.
Questions regarding the scaling limit of the standard homogenization commutator require more careful probabilistic assumptions and are addressed in the forthcoming contribution~\cite{DFG1} in the case of correlated Gaussian fields (see below for an informal discussion of these results).
In~\cite{DO1}, we further explain how this whole theory of fluctuations naturally extends to higher orders.
We refer to the introduction of the companion article~\cite{DGO1} for a general discussion of the literature on fluctuations in stochastic homogenization (a short discussion of the key pathwise structure is given at the end of this introduction).

\medskip

Let $\Aa$ be a stationary and ergodic random coefficient field on $\R^d$ that is bounded in the sense of
\begin{equation}\label{f.56}
|\Aa(x)\xi|\le|\xi|\qquad\mbox{for all}\;\xi\in\mathbb{R}^d\;\mbox{and}\;x\in\mathbb{R}^d,
\end{equation}
and satisfies the ellipticity property
\begin{equation}\label{f.40}
\int_{\R^d}\nabla u\cdot \Aa\nabla u\ge\lambda\int_{\R^d}|\nabla u|^2\qquad
\mbox{for all }u\in C^\infty_c(\R^d),
\end{equation}
for some $\lambda>0$; this notion of functional coercivity is weaker than pointwise ellipticity for systems.
Throughout the article we use scalar notation, but no iota in the proofs would change for systems under assumptions~\eqref{f.56} and \eqref{f.40}.
For all $\e>0$, we set $\Aa_\e:=\Aa(\frac \cdot \e)$,
and for a deterministic vector field $f\in C^\infty_c(\R^d)^d$ we consider the random family $(u_\e)_{\e>0}$ of unique Lax-Milgram solutions in $\R^d$ (which henceforth means the unique weak solutions in $\dot H^1(\R^d)$) of the rescaled problems
\begin{align}\label{eq:first-def-ups}
-\nabla\cdot \Aa_\e\nabla u_\e\,=\, \nabla\cdot f.
\end{align}
(The choice of considering an equation on the whole space rather than on a bounded set allows us to focus on fluctuations in the bulk, and avoid 
effects of boundary layers. The choice of taking a right-hand side (RHS) in divergence form allows to treat all dimensions at once.)
It is known since the pioneering work of Papanicolaou and Varadhan~\cite{PapaVara} and Kozlov~\cite{Kozlov-79} that, almost surely, $u_\e$ converges weakly (in $\dot H^1(\R^d)$)
as $\e \downarrow 0$ to the unique Lax-Milgram solution $\bar u$ in $\R^d$ of 
\begin{align}\label{eq:first-def-ubar-intro}
-\nabla\cdot \bar\Aa\nabla\bar u\,=\, \nabla\cdot f,
\end{align}
where $\bar\Aa$ is a deterministic and constant matrix that only depends on the law of $\Aa$.
More precisely, for any direction $e\in \R^d$, the projection $\bar\Aa e$ is the expectation of the flux of the corrector in the direction $e$,
\begin{align}\label{eq:coeff-homog}
\bar\Aa e\,=\,\expec{\Aa(\nabla\phi_e+e)},
\end{align}
where the corrector $\phi_e$ is the unique (up to a random additive constant) almost-sure solution of the corrector equation in $\R^d$,
$$
-\nabla\cdot \Aa(\nabla\phi_e+e)\,=\,0,
$$
in the class of functions the gradient of which is stationary, has finite second moment, and has zero expectation. We denote by $\phi=(\phi_i)_{i=1}^d$ the vector field the entries of which are the correctors $\phi_i$ in the canonical directions $e_i$ of $\R^d$.

\medskip

In~\cite{DGO1}, we developed a complete theory of fluctuations in stochastic homogenization for the random conductance model (see also \cite{GuM} for heuristic arguments).
The key in our theory is to focus on the so-called {\it homogenization commutator of the solution},
\begin{equation}\label{eq:homog-comm}
\Aa_\e\nabla u_\e-\bar\Aa\nabla u_\e,
\end{equation}
and to study its relation to the 
\emph{standard homogenization commutator} $\Xi:=(\Xi_i)_{i=1}^d$, where the solution $u_\e$ is replaced by $\Aa$-harmonic coordinates $x\mapsto x_i+\phi_i(x)$,
\begin{align}\label{intro:HC}
\Xi_i \,:=\,\Aa (\nabla \phi_i+\ee_i)-\bar\Aa(\nabla \phi_i+\ee_i),\qquad \Xi_{ij}:=(\Xi_i)_j.
\end{align}
In the framework of \cite{DGO1}, we showed
the following three crucial properties (which we reformulate here in the non-symmetric continuum setting):
\begin{enumerate}[(I)]
\item First and most importantly,
the two-scale expansion of the homogenization commutator of the solution
\begin{align}\label{eq:2-scale-commut}
\Aa_\e \nabla u_\e-\bar\Aa \nabla u_\e-\expec{\Aa_\e \nabla u_\e-\bar\Aa \nabla u_\e}\,\approx\,\Xi_i(\tfrac \cdot \e) \nabla_i \bar u
\end{align}
is accurate in the fluctuation scaling in the sense that for all $g\in C^\infty_c(\R^d)^d$ and $q<\infty$,
\begin{multline}\label{eq:2-scale-commut-quant}
\hspace{1.1cm}\expec{\Big| \int_{\R^d} g\cdot \big(\Aa_\e \nabla u_\e-\bar\Aa \nabla u_\e-\expec{\Aa_\e \nabla u_\e-\bar\Aa \nabla u_\e}\big)- \int_{\R^d} g\cdot \Xi_i(\tfrac \cdot \e) \nabla_i \bar u\Big|^q}^\frac1q\\
 \lesssim_{f,g,q}~\e\,\expec{\Big| \int_{\R^d} g\cdot \Xi_i(\tfrac\cdot\e) \nabla_i\bar u
 \Big|^q}^\frac1q,
\end{multline}
up to a $|\!\log\e|$ factor in the critical dimension $d=2$. This property is highly nontrivial and is due to the special form of the commutator~\eqref{eq:homog-comm}.
\smallskip\item Second, both the fluctuations of the field $\nabla u_\e$ and of the flux $\Aa_\e \nabla u_\e$ can be recovered through \emph{deterministic} projections of the fluctuations of the homogenization commutator~\eqref{eq:homog-comm}, which shows that no information is lost by passing to the commutator. More precisely, the following elementary identities are easily checked,
\begin{eqnarray}
\int_{\R^d}g\cdot (\nabla u_\e-\nabla\bar u)=- \int_{\R^d}(\bar\calP_H^*g) \cdot (\Aa_\e\nabla u_\e-\bar\Aa\nabla u_\e),\label{eq:rel-I1}\\
\int_{\R^d} g\cdot (\Aa_\e\nabla u_\e-\bar\Aa\nabla\bar u)= \int_{\R^d}(\bar\calP_L^*g)\cdot (\Aa_\e\nabla u_\e-\bar\Aa\nabla u_\e),\nonumber
\end{eqnarray}
in terms of the Helmholtz and Leray projections in $\Ld^2(\R^d)^d$,
\begin{gather}
\bar\calP_H:=\nabla(\nabla\cdot\bar\Aa \nabla)^{-1}\nabla\cdot,\qquad \bar\calP_L:=\Id-\bar\calP_H\bar\Aa,\nonumber\\
\bar\calP_H^*:=\nabla(\nabla\cdot\bar\Aa^* \nabla)^{-1}\nabla\cdot,\qquad \bar\calP_L^*:=\Id-\bar\calP_H\bar\Aa^*,\label{eq:proj-Helm-def}
\end{gather}
where $\bar\Aa^*$ denotes the transpose of $\bar\Aa$.
Similarly, the fluctuations of the field $\nabla\phi$ and of the flux $\Aa \nabla\phi$ of the corrector are also determined by those of the standard commutator $\Xi$ itself: indeed, the definition of $\Xi$ yields $-\nabla \cdot\bar \Aa \nabla \phi_i=\nabla \cdot \Xi_i$ and $\Aa (\nabla \phi_i+\ee_i)-\bar\Aa \ee_i= \Xi_i+\bar\Aa \nabla \phi_i$, to the effect of $\nabla \phi_i=-\bar\calP_H\Xi_i$ and $\Aa(\nabla \phi_i+\ee_i)-\bar\Aa\ee_i=(\Id-\bar\Aa\bar\calP_H)\Xi_i$ in the stationary sense, hence formally,
\begin{eqnarray}
 \int_{\R^d} \F:\nabla\phi(\tfrac\cdot\e)&=&- \int_{\R^d} \bar\calP_H^*\F:\Xi(\tfrac\cdot\e),\nonumber\\
 \int_{\R^d} \F:\big(\Aa_\e(\nabla \phi(\tfrac\cdot\e)+\Id)-\bar\Aa\big)&=& \int_{\R^d}\calP_L^*\F:\Xi(\tfrac\cdot\e),\label{eq:rel-J1J2}
\end{eqnarray}
where $\bar\calP_H^*$ and $\bar\calP_L^*$ act on the second index of the tensor field $\F$; a suitable sense to these identities is given in Corollary~\ref{cor:pathwise}.
\smallskip\item Third, the standard homogenization commutator $\Xi$ is an approximately local function of the coefficients $\Aa$, which allows to infer the large-scale behavior of $\Xi$ from the large-scale behavior of $\Aa$ itself. This locality property is best seen when formally computing the so-called ``vertical'' derivatives of $\Xi$ with respect to $\Aa$: Letting $\phi^*$ denote the corrector associated with the pointwise transpose field $\Aa^*$, and letting $\sigma^*$ denote the corresponding flux corrector (cf.~\eqref{f.5}), we obtain (cf.~\cite[equation~(1.10)]{DGO1} and~\eqref{eq:decomp-der-J0} below)
\begin{multline*}
\hspace{1cm}\frac{\partial}{\partial\Aa(x)}\Xi_{ij} = (\nabla \phi_j^*+\ee_j)\cdot\frac{\partial\Aa}{\partial\Aa(x)}(\nabla\phi_i+\ee_i)\\
-\nabla\cdot\bigg(\phi_j^*\frac{\partial\Aa}{\partial\Aa(x)}(\nabla\phi_i+\ee_i)\bigg)-\nabla\cdot\bigg((\phi_j^*\Aa+\sigma^*_j)\frac{\partial \nabla\phi_i}{\partial\Aa(x)}\bigg).
\end{multline*}
In view of $\frac{\partial\Aa}{\partial\Aa(x)}=\delta(\cdot-x)$, the first right-hand side term reveals an exactly local dependence upon $\Aa$. The second term is exactly local as well, but since it is written in divergence form its contribution is negligible when integrating on large scales.
The only non-local effect comes from the last term due to $\frac{\partial \nabla\phi}{\partial\Aa}$, which is given by the mixed derivative of the Green's function for $-\nabla\cdot\Aa \nabla$ and thus is expected to have only borderline integrable decay. However, it also appears inside a divergence, hence it is negligible when integrated on large scales.
\end{enumerate}
%

\smallskip

Let us comment on the structure of fluctuations revealed in~(I)--(II).
Together with the two-scale expansion~\eqref{eq:2-scale-commut-quant} of commutators, identities~\eqref{eq:rel-I1} and~\eqref{eq:rel-J1J2} imply that the fluctuations of $\nabla u_\e$, $\Aa_\e \nabla u_\e$, $\nabla\phi(\tfrac\cdot\e)$, and $\Aa_\e\nabla\phi(\tfrac\cdot\e)$ are determined to leading order by those of $\Xi(\tfrac\cdot\e)$, with error estimated in a strong norm in probability. We chose to refer to this key property as the ``pathwise'' structure of fluctuations in analogy with the language of SPDEs 
in order to emphasize that this result does not only compare probability laws of different objects (possibly constructed on different probability spaces),
but compares these objects for the same realizations of the randomness (for the same ``paths''), here in form of an error estimate at the level of stretched exponential moments.
As emphasized in~\cite{DGO1}, besides its theoretical importance, this \emph{pathwise structure} is bound to affect multi-scale computing and uncertainty quantification in an essential way. This result is indeed of the complexity-reducing type of the central results in homogenization, as it provides a description of fluctuations of a general solution by means of an off-line procedure using the standard commutator $\Xi$ in form of a two-scale expansion.
Next, in case of a weakly correlated coefficient field $\Aa$, we expect from property~(III) that $\Xi(\frac\cdot\e)$ displays the CLT scaling and that $\e^{-d/2}\Xi(\frac\cdot\e)$ converges to a white noise; the pathwise structure~(I)--(II) then allows to recover the known scaling limit results for the different quantities of interest in stochastic homogenization, as indeed shown in~\cite{DGO1} for the random conductance model.

\medskip
In the present contribution, we focus on the pathwise structure~(I)--(II). More precisely, we mainly consider the class of Gaussian coefficient fields with a covariance function that decays at infinity at some fixed (yet arbitrary) algebraic rate $(1+|x|)^{-\beta}$ parametrized by $\beta>0$,
and we show that properties~(I)--(II) still hold for this whole Gaussian class, 
which illustrates the surprising robustness of the pathwise structure with respect to the large-scale behavior of the homogenization commutator.
Indeed, in dimension $d=1$ (in which case the quantities under investigation are simpler and explicit\footnote{In dimension $d=1$, the homogenization commutator indeed simply takes the form $\Xi(x)= \bar \Aa(1 -\frac{\bar \Aa}{\Aa(x)})$, which is exactly local wrt $\Aa$.}), two typical behaviors have been identified in terms of the scaling limit of the standard homogenization commutator $\Xi$, depending on the parameter~$\beta$ (cf.~\cite{BGMP-08}),
\begin{enumerate}[\quad$\bullet$]
\item For $\beta>d=1$: The standard commutator $\Xi$ displays the CLT scaling and its rescaling $\e^{-\frac 12}\Xi(\frac\cdot\e)$ converges in law to a non-degenerate white noise (Gaussian fluctuations, local limiting covariance structure), but the convergence rate is arbitrarily slow as $\beta$ gets closer to $d=1$.
\item For $0<\beta<d=1$: The suitable rescaling $\e^{-\frac\beta2}\Xi(\frac\cdot\e)$ converges along a subsequence to a fractional Gaussian field (Gaussian fluctuations, nonlocal limiting covariance structure, potentially no uniqueness of the limit). Note that a different, non-Gaussian behavior may also occur in degenerate cases (cf.~\cite{Gu-Bal-12,LNZH-17} and second item in Remark~\ref{rem:main}).
\end{enumerate}
In particular, the pathwise result is shown to hold in both examples whereas the rescaled standard commutator does not necessarily converge to white noise or may even not converge at all.
The identification of the scaling limit of the standard commutator is thus a separate question and is addressed in~\cite{DFG1} in all dimensions for the whole range of values of $\beta>0$,
combining Malliavin calculus with techniques developed in~\cite{GNO-reg}. More precisely, this work extends~\cite{BGMP-08} to dimensions $d>1$ in the following sense,
\begin{enumerate}[\quad$\bullet$]
\item For $\beta > d$: The rescaled commutator $\e^{-\frac d2}\Xi(\frac\cdot\e)$ converges in law to a generically non-degenerate white noise.
\item For $\beta < d$: The rescaled commutator $\e^{-\frac \beta 2}\Xi(\frac\cdot\e)$ converges along a subsequence to a generically non-degenerate fractional Gaussian field. Different limits can indeed be reached in general, unless the covariance function has a self-similar profile at infinity.
\end{enumerate}
These results illustrate the fact that the standard commutator $\Xi$ is an approximately local function of the random coefficient field $\Aa$ (cf.~(III) above), which essentially allows to relate the scaling limit of the commutator with the scaling limit of the coefficient field itself (as in dimension~$d=1$).
Interestingly, this also shows that the pathwise structure of fluctuations can in general not be reduced to a quantitative joint convergence in law since there might not even be any convergence in law to talk about in the first place.

\medskip
Although we focus here for shortness on the model case of Gaussian coefficient fields, the arguments that we provide in this contribution are robust enough to cover the general setting of mutiscale functional inequalities introduced and studied in \cite{DG1,DG2}, and therefore to treat all the models of random coefficient fields considered in the reference textbook~\cite{Torquato-02} on heterogenous materials (see indeed third item of Remark~\ref{rem:main}). This makes the results of this contribution not only of theoretical but also of practical interest.

\medskip

Let us conclude this introduction with a short discussion of the recent literature concerning (I)--(III); we refer to~\cite[Section~1.4]{DGO1} for more detail. The pathwise structure~(I)--(II) of fluctuations, which we extend here to the continuum setting with long-range correlations, was first formulated and proved by us in~\cite{DGO1} for the random conductance model.
A related form of (I)--(II) was conjectured in \cite{GuM} within the variational and renormalization framework of~\cite{AS,AKM2,AKM-book}, but it has not been made rigorous yet (nor does it not appear in the textbook~\cite{AKM-book}).
A variational quantity related to the standard commutator can be first traced back to~\cite{AS}, whereas its canonical 
form~\eqref{intro:HC} used here was independently introduced in~\cite{AKM2,AKM-book} and~\cite{DGO1} (there  motivated by the seminal works of Murat and Tartar).
The locality property of the standard commutator $\Xi$ and its convergence to white noise were first established in~\cite{DGO1} for the random conductance model, and in~\cite{AKM2,GO4} for the continuum setting with a finite range of dependence assumption, while the case of long-range correlations is first considered
in our companion article \cite{DFG1}.

\section{Main results and structure of the proof}

\subsection{Notation and statement of the main results}

For some $k\ge1$ let $a$ be an $\R^k$-valued Gaussian random field, constructed on a probability space $(\Omega,\Pm)$ (with expectation~$\mathbb{E}$), which is stationary and centered, and thus characterized by its covariance function
\begin{equation*}
c(x):=\expec{a(x)\otimes a(0)},\qquad c:\R^d\to\R^{k\times k}.
\end{equation*}
We assume that the covariance function decays algebraically at infinity in the sense that there exist $\beta,C_0>0$ such that
for all $x\in \R^d$,
\begin{equation}\label{g.3}
\frac1{C_0} (1+|x|)^{-\beta} \,\le \, |c(x)| \,\le\,C_0 (1+|x|)^{-\beta}.
\end{equation}
Given a map $h\in C^1_b(\R^k)^{d\times d}$, we define $\Aa:\R^d\to\R^{d\times d}$ by $\Aa(x)=h(a(x))$, and assume that it satisfies the conditions~\eqref{f.56} and~\eqref{f.40} almost surely. We then (abusively) call the coefficient field $\Aa$ \emph{Gaussian with parameter $\beta>0$}.
If $\Aa$ is Gaussian with parameter $\beta$, then $\Aa$ is ergodic, hence we have existence and uniqueness of correctors $\phi$ and of the homogenized coefficients $\bar\Aa$ (cf.~Lemma~\ref{si} below).
From a technical point of view, we shall rely on (and frequently refer to) results and methods developed in~\cite{DG1,GNO-reg,GNO-quant}.

\medskip
Throughout the article, we use the notation $\lesssim_{(\dots)}$ (resp. $\gtrsim_{(\dots)}$) for $\le C\times$ (resp. $\ge C\times$), where the multiplicative constant
$C$ depends on $d,\lambda,\beta,\|\nabla h\|_{\Ld^\infty}$, on the constant $C_0$ in~\eqref{g.3}, and on the additional parameters ``$(\dots)$'' if any. We write $\simeq_{(\dots)}$ when both $\lesssim_{(\dots)}$ and $\gtrsim_{(\dots)}$ hold. In an assumption, we use the notation $\ll_{(\dots)}$ for $\le\frac1C\times$ for some (large enough) constant $C\simeq_{(\dots)}1$.

\medskip
We now define a string of random functionals that encode the fluctuations of the different objects of interest.
The notation $I$ is reserved to functionals involving the solution operator,
and the notation $J$ to functionals involving correctors; the subscript $_0$ is reserved to commutators,
the subscript $_1$ to fields, and the subscript $_2$ to fluxes.
We consider the fluctuations of the commutator $\Aa_\e \nabla u_\e - \bar \Aa \nabla u_\e$, of the field $\nabla u_\e$, and of the flux $\Aa_\e\nabla u_\e$ of the solution to~\eqref{eq:first-def-ups}, as encoded by the (centered) random bilinear functionals $I_0^\e:(f,g)\mapsto I_0^\e(f,g)$, $I_1^\e:(f,g)\mapsto I_1^\e(f,g)$, and $ I_2^\e:(f,g)\mapsto I_2^\e(f,g)$ defined
for all $f,g\in C^\infty_c(\R^d)^d$ by
\begingroup\allowdisplaybreaks
\begin{eqnarray*}
I_0^\e(f,g)&:=& \int_{\R^d}g \cdot (\Aa_\e \nabla u_\e - \bar \Aa \nabla u_\e-\expec{\Aa_\e \nabla u_\e-\bar\Aa \nabla u_\e}),\\
 I_1^\e(f,g)&:=& \int_{\R^d}g\cdot\nabla(u_\e-\expec{u_\e}),\\
 I_2^\e(f,g)&:=& \int_{\R^d}g\cdot\big(\Aa_\e\nabla u_\e-\expec{\Aa_\e\nabla u_\e}\big).
\end{eqnarray*}
Likewise, we consider the fluctuations of the standard commutator $\Xi=(\Aa -\bar \Aa)(\nabla \phi+\Id)$, of the corrector field $\nabla\phi$, and of the corrector flux $\Aa(\nabla\phi+\Id)$ as encoded by the (centered) random linear functionals $J_0^\e:\F\mapsto J_0^\e(\F)$, $J_1^\e:\F\mapsto J_1^\e(\F)$, and $ J_2^\e:\F\mapsto J_2^\e(\F)$ defined for all $\F\in C^\infty_c(\R^d)^{d\times d}$ by
\begin{eqnarray*}
 J_0^\e(\F)&:=& \int_{\R^d} \F(x): \Xi(\tfrac x\e)\,dx,\\
 J_1^\e(\F)&:=& \int_{\R^d} \F(x): \nabla \phi(\tfrac{x}{\e})\,dx,\\
 J_2^\e(\F)&:=& \int_{\R^d} \F(x):\big(\Aa_\e(x) (\nabla \phi(\tfrac{x}{\e})+\Id)-\bar\Aa\big)\,dx.
\end{eqnarray*}
\endgroup

We first prove the following boundedness result for $ J_0^\e$, establishing the suitable $\beta$-dependent scaling for the fluctuations of the homogenization commutator
(see also~\cite[Theorem~1]{GNO-quant}).
More precisely, in the spirit of~(III), this shows that large-scale averages of the standard commutator have the same scaling $\pi_*(\tfrac 1\e)^{-1/2}$ as large-scale averages of the coefficient field $\Aa$ itself (cf.~\cite[Proposition~1.5]{DG1}); in the case of integrable correlations, this is the CLT scaling~$\e^{d/2}$.

\begin{prop1}[Fluctuation scaling]\label{prop:scaling}
Let $d\ge1$, assume that the coefficient field $\Aa$ is Gaussian with parameter $\beta>0$, define $\pi_*:\R^+\to\R^+$ by
\begin{align}\label{e.pi*}
\pi_*(t)\,:=\,
\left\{\begin{array}{lll}
(1+t)^\beta&:&\beta<d,\\
(1+t)^d\, \frac1{\log(2+t)}&:&\beta=d,\\
(1+t)^d&:&\beta>d ,
\end{array}\right.
\end{align}
and define the rescaled functional
\[\widehat J_{0}^\e\,:=\,\pi_*(\tfrac1\e)^\frac12 J_{0}^\e.\]
For all $0<\e\le1$, $\F\in C^\infty_c(\R^d)^{d\times d}$, $0<p-1\ll1$, and $\alpha>\frac{d-(\beta\wedge d)}2+d\frac{p-1}{2p}$, we have
\begin{equation}\label{t1.var}
|\widehat J_{0}^\e(\F)|\,\le\, \calC_{\alpha,p}^{\e,\F}\,\big(\|w_1^{\alpha}F\|_{\Ld^{2p}}+\mathds1_{\beta\le d}\|[F]_2\|_{\Ld^p}\big),
\end{equation}
where $w_1(z):=1+|z|$, $[F]_2(x):=(\fint_{B(x)}|F|^2)^\frac12$, and where $\calC_{\alpha,p}^{\e,\F}$ is a random variable with stretched exponential moments: there exists $\gamma_1\simeq1$ such that
\[\sup_{0<\e<1}\expec{\exp\Big(\frac{1}{C_{\alpha,p}}(\calC_{\alpha,p}^{\e,\F})^{\gamma_1}\Big)}\,\le\,2\]
for some (deterministic) constant $C_{\alpha,p}\simeq_{\alpha,p}1$.
\end{prop1}

Our next main result establishes the accuracy of the two-scale expansion error for large-scale averages of the homogenization commutator in the suitable fluctuation scaling. This error is encoded by the following (centered) random bilinear functional,
\begin{multline}\label{e.2sError}
E^\e(f,g)\,:=\, \int_{\R^d}g \cdot \big(\Aa_\e \nabla u_\e - \bar \Aa \nabla u_\e-\expec{\Aa_\e \nabla u_\e-\bar\Aa \nabla u_\e}\big)-
 \int_{\R^d}g \cdot \Xi(\tfrac \cdot \e) \nabla \bar u\\
=\,I_0^\e(f,g)-J_0^\e( \nabla \bar u \otimes g)
\,=\,I_0^\e(f,g)-J_0^\e(\bar \calP_H f \otimes g).
\end{multline}
More precisely, we show that the typical scaling of this error $I_0^\e(f,g)-J_0^\e(\bar \calP_H f \otimes g)$ is an order $\e\mu_*(\frac1\e)$ smaller than the typical scaling of large-scale averages of the commutator $J_0^\e(\bar \calP_H f \otimes g)$ itself.
In view of the generic non-degeneracy result in~\cite{DFG1}, this can be viewed as a relative error estimate.
This property summarizes the pathwise structure of fluctuations and is the key part of our theory.
\begin{theor}[Pathwise structure of fluctuations]\label{th:main}
Let {$d\ge1$}, assume that the coefficient field $\Aa$ is Gaussian with parameter $\beta>0$, let $\pi_*$ be defined in~\eqref{e.pi*}, and define $\mu_*:\R_+\to\R_+$ by
\begin{equation}\label{mudbeta}
\mu_*(r)\,:=\, 
\left\{
\begin{array}{lll}
1 &:&\beta>2,\,d>2,\\
\log^\frac12 (2+r)&:&\beta>2,\,d=2,\\
\sqrt{1+r} &:&\beta>1,\,d=1,\\
\log^\frac12(2+r)&:&\beta=2,\,d> 2,\\
\log(2+r)&:&\beta=2,\,d= 2,\\
\sqrt{1+r}\,\log^\frac12 (2+r)&:&\beta=1,\,d=1,\\
(1+r)^{1-\frac\beta2}&:&\beta<2,\,d\ge 2,~\text{or}~\beta<1,\,d=1.
\end{array}
\right.
\end{equation}
Set $\mu_*(z):=\mu_*(|z|)$, recall the notation $w_1(z):=1+|z|$, and consider the rescaled error functional $\widehat E^\e:=\pi_*(\tfrac1\e)^\frac12 E^\e$.
For all $0<\e\le1$, $f,g\in C^\infty_c(\R^d)^{d}$, $0<p-1\ll1$, and~$\alpha>\frac{d-(\beta\wedge d)}2+d\frac{p-1}{4p}$, we have
\begin{multline}\label{t1.path}
|\widehat E^\e(f,g)|
\,\le \, \e \mu_{*}(\tfrac1\e)\,\calC^{\e,f,g}_{\alpha,p}\Big(\|\mu_*\nabla f\|_{\Ld^4}\|w_1^{\alpha}g\|_{\Ld^{4p}}+\|\mu_*\nabla g\|_{\Ld^4}\|w_1^{\alpha}f\|_{\Ld^{4p}}\\
+\mathds1_{\beta\le d}\big(\|\mu_*\nabla f\|_{\Ld^2}\|g\|_{\Ld^2\cap\Ld^{2p}}+\|\mu_*\nabla g\|_{\Ld^2}\|f\|_{\Ld^2\cap\Ld^{2p}}\big)\Big),
\end{multline}
where $\calC_{\alpha,p}^{\e,f,g}$ is a random variable with stretched exponential moments: there exists $\gamma_2\simeq1$
such that
\[\sup_{0<\e<1}\expec{\exp\Big(\frac{1}{C_{\alpha,p}}(\calC_{\alpha,p}^{\e,f,g})^{\gamma_2}\Big)}\,\le\,2\]
for some (deterministic) constant $C_{\alpha,p}\simeq_{\alpha,p}1$.
\end{theor}
\begin{rem}\label{rem:main}\mbox{}
\begin{enumerate}[\quad$\bullet$]
\item The exponents $\gamma_1$ and $\gamma_2$ in the above results can be made explicit; we do not pursue this direction since
the values obtained in the proofs are not expected to be optimal.
\smallskip
\item The $\e$-scaling in the above results is believed to be optimal.
The rescaling in the definition of $\widehat J_0^\e$ and $\widehat E^\e$ is natural since it precisely coincides with the scaling of large-scale averages of the coefficient field $\Aa$ itself.
For some non-generic examples, the bound~\eqref{t1.var} may however overestimate the variance. In dimension $d=1$, one may indeed construct explicit Gaussian coefficient fields $\Aa$ such that fluctuations of the homogenization commutator $J_0^\e$ are of smaller order than what~\eqref{t1.var} predicts~\cite{Taqqu,Gu-Bal-12,LNZH-17}, in which case the suitable rescaling of $J_0^\e$ has a non-Gaussian limit.
In such situations, the pathwise property~\eqref{t1.path} (or its higher-order pathwise version as in \cite{DO1}) might still provide relevant information. General necessary and sufficient conditions for the sharpness of~\eqref{t1.var} are provided in~\cite{DFG1}.
\smallskip
\item The proofs of the above results are robust enough to cover the general setting of multiscale functional inequalities introduced in~\cite{DG1,DG2}.
In the case of functional inequalities with oscillation, we may indeed use Cauchy-Schwarz' inequality and an energy estimate to replace the perturbed functions $\tilde\phi$ and $\nabla\tilde u$ appearing in the representation formula~\eqref{eq:corr-eqn-der-pre} below by their unperturbed versions $\phi$ and $\nabla u$. 
This allows to conclude whenever the weight has a superalgebraic decay (see indeed~\cite[proof of Theorem~4]{GNO-quant}). 
If one is only interested in Gaussian coefficient fields, one may replace the use of functional inequalities by a direct use of the Brascamp-Lieb inequality in terms of Malliavin calculus, which allows
to shorten some of the proofs (and improve the norms of the test functions $F,f,g$), cf.~\cite{DO1}.
\qedhere
\end{enumerate}
\end{rem}
In view of the identities~\eqref{eq:rel-I1} and~\eqref{eq:rel-J1J2}, the above pathwise result implies that the large-scale fluctuations of $I_0^\e$, $I_1^\e$, $I_2^\e$, $J_1^\e$, and $J_2^\e$ are driven by the fluctuations of $J_0^\e$ in a pathwise sense (see~\cite[Corollary~2.4]{DGO1} for details).

\begin{cor}[\cite{DGO1}]\label{cor:pathwise}
Let $d\ge2$, assume that the coefficient field $\Aa$ is Gaussian with parameter $\beta>0$, let $\pi_*$ and $\mu_*$ be defined by~\eqref{e.pi*} and~\eqref{mudbeta}, let $\bar\calP_H$, $\bar\calP_H^*$, and $\bar\calP_L^*$ be as in~\eqref{eq:proj-Helm-def}, and recall the rescaled functionals
\[\widehat I_i^\e:=\pi_*(\tfrac1\e)^\frac12 I_i^\e,\qquad\widehat J_i^\e:=\pi_*(\tfrac1\e)^\frac12J_i^\e,\qquad i=0,1,2.\]
For all $\e>0$ and $f,g\in C^\infty_c(\R^d)^{d}$, we have for all $0<p-1\ll1$ and~$\alpha>\frac{d-(\beta\wedge d)}2+d\frac{p-1}{4p}$,
\begin{multline*}
|\widehat I_0^\e(f,g)-\widehat J_0^\e(\bar\calP_Hf\otimes g)|+|\widehat I_1^\e(f,g)-\widehat J_0^\e(\bar\calP_Hf\otimes \bar\calP_H^*g)|+|\widehat I_2^\e(f,g)+\widehat J_0^\e(\bar\calP_Hf\otimes\bar\calP_L^*g)|\\
\,\le \, \e \mu_{*}(\tfrac1\e)\,\calC^{\e,f,g}_{\alpha,p}\Big(\|\mu_*\nabla f\|_{\Ld^4}\|w_1^{\alpha}g\|_{\Ld^{4p}}+\|\mu_*\nabla g\|_{\Ld^4}\|w_1^{\alpha}f\|_{\Ld^{4p}}\\
+\mathds1_{\beta\le d}\big(\|\mu_*\nabla f\|_{\Ld^2}\|g\|_{\Ld^2\cap\Ld^{2p}}+\|\mu_*\nabla g\|_{\Ld^2}\|f\|_{\Ld^2\cap\Ld^{2p}}\big)\Big),
\end{multline*}
where $\calC_{\alpha,p}^{\e,f,g}$ is a random variable with stretched exponential moments independent of $\e$ as in the statement of Theorem~\ref{th:main}.
In addition, for all $\e>0$ and $\F\in C^\infty_c(\R^d)^{d\times d}$, we have almost surely
\begin{align*}
J_1^\e(\F)=-J_{0}^\e(\bar\calP_H^*\F),\qquad J_2^\e(\F)=J_{0}^\e(\bar\calP_L^*\F),
\end{align*}
where in particular we may give
an almost sure meaning to $J_0^\e(\bar\calP_H^*\F)$ and $J_{0}^\e(\bar\calP_L^*\F)$ for all $\F\in C^\infty_c(\R^d)^{d\times d}$, even when $\bar\calP_H^*\F$ and $\bar\calP_L^*\F$ do not have integrable decay.
\end{cor}


\subsection{Structure of the proof}

We describe the string of arguments that leads to Proposition~\ref{prop:scaling} and Theorem~\ref{th:main}.
Next to the corrector $\phi$, we first need to recall the notion of flux corrector $\sigma$, which was recently introduced in the stochastic setting in~\cite[Lemma~1]{GNO-reg}
and allows to put the equation for the two-scale homogenization error in divergence form (cf.~\eqref{eq:wf} below).
The extended corrector $(\phi,\sigma)$ is only defined up to an additive constant, and we henceforth choose the anchoring $\fint_B(\phi,\sigma)=0$ on the centered unit ball~$B$.
\begin{lem}[\cite{GNO-reg}]\label{si}
Let the coefficient field $\Aa$ be stationary and ergodic. Then there exist two random tensor fields
$(\phi_i)_{1\le i\le d}$ and $(\sigma_{ijk})_{1\le i,j,k\le d}$ with the following properties:
The gradient fields $\nabla\phi_i$ and $\nabla\sigma_{ijk}$ are stationary\footnote{That is, $\nabla\phi_i(\Aa;\cdot+z)=\nabla\phi_i(\Aa(\cdot+z);\cdot)$ and $\nabla\sigma_{ijk}(\Aa;\cdot+z)=\nabla\sigma_{ijk}(\Aa(\cdot+z);\cdot)$ a.e.\@ in $\R^d$, for all shift vectors $z\in\mathbb{R}^d$.}
and have finite second moments and vanishing expectations:
\begin{equation}\label{si.2}
\expec{|\nabla\phi_i|^2}\le \frac{1}{\lambda^2},\quad
\sum_{j,k=1}^d\expec{|\nabla\sigma_{ijk}|^2}
\le {4d}\Big(\frac{1}{\lambda^2}+1\Big),\quad
\expec{\nabla\phi_i}=\expec{\nabla\sigma_{ijk}}=0.
\end{equation}
Moreover, for all $i$, the field $\sigma_i:=(\sigma_{ijk})_{1\le j,k\le d}$ is skew-symmetric, that is,
\begin{equation}\label{f.19}
\sigma_{ijk}=-\sigma_{ikj}.
\end{equation}
Finally, the following equations are satisfied a.s.\@ in the distributional sense on $\R^d$,
\begin{eqnarray}
-\nabla\cdot \Aa(\nabla\phi_i+e_i)&=&0,\label{f.2}\\
\nabla\cdot\sigma_i&=&q_i-\expec{q_i},\label{f.5}\\
-\triangle\sigma_{ijk}&=&\partial_jq_{ik}-\partial_kq_{ij},\nonumber
\end{eqnarray}
where $q_i=(q_{ij})_{1\le j\le d}$ is given by $q_i:=\Aa(\nabla\phi_i+e_i)$,
and where the (distributional) divergence of a tensor field is defined as $(\nabla\cdot\sigma_i)_j:=\sum_{k=1}^d\nabla_k\sigma_{ijk}$.
\end{lem}

\medskip

The proofs of Proposition~\ref{prop:scaling} and Theorem~\ref{th:main} are based on the combination of 
three main ingredients:
\begin{itemize}
\item A sensitivity calculus combined with functional inequalities for Gaussian ensembles~\cite{DG1,DG2};
\item The bounds on correctors proved in \cite{GNO-quant};
\item A duality argument combined with the large-scale
(weighted) Calder\'on-Zygmund estimates of \cite{GNO-reg}.
\end{itemize}
In the case when the coefficients satisfy a finite range of dependence assumption rather than a functional inequality, 
we do not have a convenient sensitivity calculus at our disposal, and this first ingredient can be replaced
by a semi-group approach that provides a convenient disintegration of scales, cf.~\cite{DFG2}.

\medskip

The sensitivity calculus measures the influence of changes of the coefficient field $\Aa$ 
on random variables $X=X(\Aa)$ via the functional (or Malliavin-type) derivative $\partial^\fun X(x)=\frac{\partial X}{\partial \Aa}(\Aa,x)$,
that is, the $\Ld^2(\mathbb{R}^d)^{d\times d}$-gradient of $X$ wrt $\Aa$.
We recall that this functional derivative is characterized as follows, for any compactly supported perturbation $\bb\in\Ld^\infty(\R^d)^{d\times d}$,
\begin{equation}
\int_{\mathbb{R}^d}\partial^\fun X(\Aa,x): \bb(x)\,dx\,:=\,\lim_{t\downarrow 0}\frac{1}{t}\big(X(\Aa+t\bb)-X(\Aa)\big). \label{s.50}
\end{equation}
This quantity measures the sensitivity of the random variable $X=X(\Aa)$ wrt changes in the coefficient
field.
This sensitivity calculus is a building block to control the variance and the entropy of $X$ via functional inequalities in the probability space~\cite{DG1}.
A crucial role is played by the parameter $\beta>0$ that characterizes the decay of the covariance function of $\Aa$,
and we define as follows a weighted norm $\Norm{\cdot}^2_\beta$ on random fields $G$, depending on $\beta>0$,
\begin{eqnarray}
\Norm{G}^2_\beta &:=& \int_1^\infty \| G\|_{\ell}^2\,\,\ell^{-\beta-1}\,d\ell,\label{e.operator-norm}
\end{eqnarray}
where for all $\ell\ge 1$
\begin{eqnarray}
 \| G\|_{\ell}^2&:=&\ell^{-d}\int_{\R^d}\Big( \int_{B_{\ell}(z)} |G|\Big)^2 dz.\label{eq:partition:0}
 \end{eqnarray}
As shown in~\cite[Proposition~2.4]{DG1}, in the integrable case $\beta>d$, we can drop the integral over $\ell$,
in which case
\begin{align}
\Norm{G}^2_\beta~\simeq~ \Norm{G}^2 ~:= \|G\|_1^2.\label{e.operator-norm-simpl}
\end{align}
In these terms, we may formulate the following multiscale logarithmic Sobolev inequality for the Gaussian coefficient field $\Aa$.
In view of~\eqref{e.operator-norm-simpl}, for $\beta>d$, this reduces to the standard logarithmic Sobolev inequality (LSI).
The proof is based on a corresponding Brascamp-Lieb inequality (cf.~\cite[Theorem~3.1]{DG2}).
\begin{lem}[\cite{DG2}]\label{lem:LSI}
Assume that the coefficient field $\Aa$ is Gaussian with parameter $\beta>0$. Then
for all random variables $X=X(\Aa)$,
\[\ent{X^2}\,:=\,\expec{X^2 \log X^2}-\expec{X^2}\expec{\log X^2}
\,\lesssim\,\expec{\Norm{\partial^\fun X}^2_\beta}.\qedhere\]
\end{lem}

\medskip

Our general strategy for the proof of Proposition~\ref{prop:scaling} and Theorem~\ref{th:main} consists in estimating the weighted norm~\eqref{e.operator-norm} of the functional derivatives of $J_0^\e(\F)$ and of $E^\e(f,g)$.
The following lemma provides a useful representation formula for these functional derivatives. This is a continuum version of~\cite[Lemma~3.2]{DGO1}. By scaling, it is enough to consider $\e=1$, and we write for simplicity $J_0:=J_0^1$ and $E:=E^1$.
\begin{lem}[Representation formulas]\label{lem:decompI3eps}
Let the coefficient field $\Aa$ be Gaussian with parameter $\beta>0$.
For all $f\in C^\infty_c(\R^d)^d$, let $u:=u_1$ denote the solution of~\eqref{eq:first-def-ups} (with $\e=1$), let $\bar u$ denote the solution of~\eqref{eq:first-def-ubar-intro},
and define the two-scale expansion error $w_{f}:=u-(1+\phi_i\nabla_i)\bar u$.
Then, for all $\F\in C^\infty_c(\R^d)^{d\times d}$,
\begin{align}\label{eq:decomp-der-J0}
\partial^\fun J_0(\F)=(\F_{ij}\,\ee_j+\nabla S_i)\otimes(\nabla\phi_i+\ee_i),
\end{align}
and for all $g\in C^\infty_c(\R^d)^d$,
\begin{multline}\label{eq:decomp-der-E}
\partial^\fun E(f,g)=g_{j}\,(\nabla\phi_j^*+\ee_j)\otimes(\nabla w_{f}+\phi_i\nabla\nabla_i\bar u)
+(\phi_j^*\,\nabla g_{j}+\nabla r)\otimes\nabla \bar u\\
-\big(\phi_j^*\,\nabla(g_{j}\nabla_i\bar u)+\nabla R_{i}\big)\otimes(\nabla\phi_i+\ee_i),
\end{multline}
where the auxiliary fields $S=(S_i)_{i=1}^d$, $r$, and $R=(R_{i})_{i=1}^d$ are the Lax-Milgram solutions in~$\R^d$ of
\begin{eqnarray}
-\nabla\cdot\Aa^*\nabla S_i&=&\nabla\cdot\big(\F_{ij}(\Aa^*-\bar\Aa^*)\ee_j\big),\label{eq:aux-seps}\\
-\nabla\cdot\Aa^*\nabla r&=&\nabla\cdot\big((\phi_j^*\Aa^*-\sigma_{j}^*)\nabla g_{j}\big),\label{eq:aux-reps}\\
-\nabla\cdot\Aa^*\nabla R_{i}&=&\nabla\cdot\big((\phi_j^*\Aa^*-\sigma_{j}^*)\nabla(g_{j}\nabla_i\bar u)\big),\label{eq:aux-Reps}
\end{eqnarray}
and $\Aa^*$ denotes the pointwise transpose coefficient field of $\Aa$, and $(\phi^*,\sigma^*)$ denotes the corresponding extended corrector
(recall that $\overline{\Aa^*}=\bar\Aa^*$).
\qedhere
\end{lem}
Before we turn to the (technical) estimates of $\Norm{\partial^\fun J_0(\F)}_\beta$ and $\Norm{\partial^\fun E(f,g)}_\beta$,
let us give an informal discussion of the scalings of the terms appearing in~\eqref{eq:decomp-der-J0} and \eqref{eq:decomp-der-E}.
To keep this discussion short, assume that $\nabla \phi, \nabla \sigma$ are bounded (which only holds after taking stochastic moments), 
that $|\phi(x)|+|\sigma(x)|\lesssim \mu_*(|x|)$ (which again only holds after taking stochastic moments), and that the Helmholtz projections associated with $-\nabla\cdot\Aa^*\nabla$ (and used to define $S_i$, $r$, and $R_i$ via 
\eqref{eq:aux-seps}--\eqref{eq:aux-Reps}) enjoy perfectly local bounds in the sense that 
$$
-\nabla\cdot\tilde \Aa\nabla z = \nabla \cdot Z \qquad\implies \qquad |\nabla z(x)|\,\lesssim\, |Z(x)|\quad \text{ for all }x\in \R^d,
$$
with $\tilde \Aa=\Aa$ or $\bar \Aa$
(which even in the homogeneous case $\tilde\Aa=\bar\Aa$ would only hold after taking suitable Lebesgue norms in view of the Calderón-Zygmund theory). 
For $J_0(F)$, equation~\eqref{eq:aux-seps} would then yield the pointwise bound
$
|\partial^\fun J_0(F)| \,\lesssim\, |F|,
$
hence
\begin{eqnarray*}
\Norm{\partial^\fun J_0(F)}^2_\beta&\lesssim & \int_1^\infty\ell^{-d}\int_{\R^d}\Big( \int_{B_{\ell}(z)} |F|\Big)^2 dz\,\,\ell^{-\beta-1}\,d\ell.
\end{eqnarray*}
To estimate the RHS, assume that $F$ is compactly supported in $B_R$ for some $R>0$,
so that
$$
\ell^{-d} \int_{\R^d}\Big( \int_{B_{\ell}(z)} |F|\Big)^2 dz\,\lesssim\, (\ell^{d} \mathds{1}_{\ell \le R} +R^d \mathds{1}_{\ell>R})\Big( \int_{\R^d}|F|^2\Big),
$$
which after integration yields, in view of~\eqref{e.pi*},
\begin{equation*}
\Norm{\partial^\fun J_0(F)}^2_\beta \,\lesssim\,\Big(\int_1^R \ell^{d-\beta-1}d\ell+R^d\int_R^\infty \ell^{-\beta-1}d\ell\Big)\Big( \int_{\R^d}|F|^2\Big)\,\lesssim\,R^d \pi_*^{-1}(R) \int_{\R^d}|F|^2.
\end{equation*}
Replacing $F$ by $\e^d F(\frac \cdot \e)$, hence replacing $R$ by $\frac1\e R$, we conclude by LSI,
$$\var{\pi_*(\tfrac 1\e)^\frac12J_0^\e(\F)} \,\lesssim_R\, \|\F\|_{\Ld^2(\R^d)}^2,$$
as claimed in Proposition~\ref{prop:scaling} (with a slightly stronger norm of the test function $F$).

\medskip
We now turn to the two-scale expansion error, for which \eqref{eq:aux-reps} and \eqref{eq:aux-Reps} would yield the following pointwise bound, under the simplifying assumptions,
$$
|\partial^\fun E(f,g)|\,\lesssim\,|g|(|\nabla w_{f}|+\mu_*|\nabla^2\bar u|)
+\mu_*|\nabla g||\nabla\bar u|.
$$
Let us further reformulate the RHS. 
On the one hand, since $-\nabla \cdot \bar \Aa \nabla \bar u = \nabla \cdot f$, the simplifying assumptions
yield the pointwise bounds $|\nabla \bar u|\,\lesssim\, |f|$ and $|\nabla^2 \bar u|\lesssim |\nabla f|$.
On the other hand, the function $w_f$ satisfies the equation
$
-\nabla \cdot \Aa \nabla w_f\,=\, \nabla \cdot ((\Aa \phi-\sigma)\nabla^2 \bar u),
$
(cf.~\cite[Remark~3]{GNO-quant})
so that our simplifying assumptions yield this time $|\nabla w_f|\, \lesssim \,\mu_* |\nabla^2 \bar u|\,\lesssim\, \mu_* |\nabla f|$.
This leads to the bound
\[|\partial^\fun E(f,g)|\,\lesssim\ \mu_*(|g||\nabla f|+|\nabla g|| f|).\]
As above, after rescaling, and using $\mu_*(|\frac x\e|) \lesssim \mu_*(|x|)\mu_*(\frac1\e)$, we conclude by LSI,
\[\var{\pi_*(\tfrac1\e)^\frac12E^\e(f,g)} \,\lesssim\, \e^2 \mu_*^2(\tfrac 1\e)\big( \|\mu_* \nabla g\|_{\Ld^4(\R^d)}^2 \|f\|_{\Ld^4(\R^d)}^2 +\|\mu_* \nabla f\|_{\Ld^4(\R^d)}^2\|g\|_{\Ld^4(\R^d)}^2\big),\]
as claimed in Theorem~\ref{th:main} (with slightly stronger norms of the test functions).
Note that the additional factor $\e^2$ comes from the gradients $\nabla f$ and $\nabla g$ in the bound 
of the functional derivative (which indeed both yield an $\e$ factor by rescaling).

\medskip

The rigorous proof of Proposition~\ref{prop:scaling} and Theorem~\ref{th:main} 
amounts to taking care of the fact that the above simplifying assumptions only hold in weaker forms.
More precisely:
\begin{itemize}
\item Bounds on the correctors only hold for stretched exponential moments and not pointwise
(cf.~Lemma~\ref{lem:bd-sigma} below), so that the bounds on $|\partial^\fun J_0(\F)$ and $|\partial^\fun E(f,g)|$ will not hold pointwise but for stretched exponential moments.
\smallskip\item More importantly, the Helmholtz projection never enjoys pointwise bounds, which must be weakened in two ways.
First, for the homogeneous operator $-\nabla\cdot\bar\Aa\nabla$, we must resort to the boundedness of the Helmholtz projection in $\Ld^p$ spaces for $1<p<\infty$ (Calder\'on-Zygmund estimates). Second, for the heterogeneous operator $-\nabla\cdot\Aa\nabla$, regularity theory can only hold on large scales~\cite{AKM-book,GNO-reg}, so that Calder\'on-Zygmund estimates must be locally averaged at some random scale $r_*$ (cf.~Lemma~\ref{lem:r*}); we will have to get rid of this random local average at some point using H\"older's inequality and a small weight (see e.g.~the second RHS factor in~\eqref{est1-LSI}). Finally, when the corrector does not have uniformly bounded moments (that is, when it grows at infinity), we further need to resort to weighted Calder\'on-Zygmund estimates (cf.~Lemma~\ref{lem:r*}(c)); see e.g.~the weight $\mu_*$ in the third RHS factor in~\eqref{est2-LSI}.
\end{itemize}

\medskip
Estimates on $\Norm{\partial^\fun J_0(\F)}_\beta$ and $\Norm{\partial^\fun E(f,g)}_\beta$ are obtained in the following two technical propositions.
As above, we prove the estimate for $\e=1$ and then argue by scaling.
Since we need some flexibility in the weights, some estimates involve a parameter $R\ge 1$.
This parameter is arbitrary and should be thought of as being $R=\frac1\e$ for the proof of the main results (similarly as in the above informal discussion). Henceforth we write~$\int$~instead of $\int_{\R^d}$ for simplicity.

\begin{prop}[Main estimates]\label{prop:boundG2}
Let the coefficient field $\Aa$ be Gaussian with parameter $\beta>0$.
Let $\pi_*$ and $\mu_*$ be defined by~\eqref{e.pi*} and~\eqref{mudbeta}, respectively, and let the random field $r_*$ be the minimal radius of Lemma~\ref{lem:r*} below.
For $F\in C^\infty_c(\R^d)$ we denote by $[F]_2(x)\,:=\,(\fint_{B(x)} |F|^2)^\frac12$ the moving local quadratic average,
and for $R\ge1$ we set $w_R(x):=\frac{|x|}R+1$. Then the following hold:
\begin{enumerate}[(i)]
\item If $\beta>d$, we have for all $R\ge1$, $0<\alpha-d\ll1$, and $0<p-1\ll_\alpha1$,
\begin{align}\label{est1-LSI}
\qquad\Norm{\partial^\fun J_0(\F)}^2\,\lesssim_{\alpha,p}\,
r_*(0)^{\alpha\frac{p-1}p}\Big(\int r_*^{d\frac p{p-1}} w_R^{-\alpha}\Big)^\frac{p-1}p\Big( \int w_R^{\alpha(p-1)}|F|^{2p}\Big)^\frac1p.
\end{align}
\item If $\beta\le d$, we have for all $R\ge1$, $0< \gamma<\beta$, $0<\alpha-d\ll1$, and $0<p-1\ll_{\gamma,\alpha}1$,
\begin{multline}\label{est1-WLSI}
\quad~\Norm{\partial^\fun J_0(\F)}_\beta^2\,\lesssim_{\alpha,p}\, R^d\pi_*(R)^{-1}\mathrm{[RHS\eqref{est1-LSI}]}
+R^{2d-\beta-\frac{2d}p}r_*(0)^{2d\frac{p-1}p}\Big( \int r_*^{d(1-\frac p2)}[F]_2^p\Big)^\frac2p\\
+R^{d-\beta}r_*(0)^{d-\gamma+\alpha\frac{p-1}p}\Big( \int r_*^{\gamma\frac{p}{p-1}}w_R^{-\alpha}\Big)^\frac{p-1}p\Big(\int w_{R}^{(d-\gamma)p+\alpha(p-1)} |\F|^{2p}\Big)^\frac1p,
\end{multline}
where we use the short-hand notation $\mathrm{[RHS\eqref{est1-LSI}]}$ for the RHS of~\eqref{est1-LSI}.\qedhere
\end{enumerate}
\end{prop}
\begin{prop}[Main estimates --- cont'd]\label{prop:boundG1}
Let the coefficient field $\Aa$ be Gaussian with parameter $\beta>0$.
Let $\pi_*$ and $\mu_*$ be defined by~\eqref{e.pi*} and~\eqref{mudbeta}, respectively, and let the random fields $r_*$ and $\calC$ be defined in Lemmas~\ref{lem:r*} and~\ref{lem:bd-sigma} below.
For $F\in C^\infty_c(\R^d)^d$ we denote by $[F]_\infty(x)\,:=\,\sup_{B(x)}|F|$ the moving local supremum,
and we recall the notation $w_R(x)=\frac{|x|}{R}+1$. Then the following hold:
\begin{enumerate}[(i)]
\item If $\beta>d$, we have for all $R\ge1$, $0<\alpha-d\ll1$, and $0<p-1\ll_\alpha1$,
\begingroup\allowdisplaybreaks
\begin{eqnarray}
\lefteqn{\Norm{\partial^\fun E(f,g)}^2\,\lesssim_{\alpha,p}\,r_*(0)^{\alpha\frac{p-1}{2p}}\Big( \int r_*^{2d\frac{2p}{p-1}}w_R^{-\alpha}\Big)^\frac{p-1}{2p}}\nonumber\\
&&\quad\times\bigg(\Big( \int \calC^{4}\mu_*^{4}[\nabla^2\bar u]_\infty^{4}\Big)^\frac{1}{2}\Big( \int w_R^{\alpha(p-1)}[g]_\infty^{4p}\Big)^\frac1{2p}\nonumber\\
&&\qquad\qquad+\,\Big( \int \calC^4\mu_*^4[\nabla g]_\infty^4\Big)^\frac12\Big( \int w_R^{\alpha(p-1)}(|f|+[\nabla\bar u]_\infty)^{4p}\Big)^\frac1{2p}\bigg).\label{est2-LSI}
\end{eqnarray}
\endgroup
\item If $\beta\le d$, we have for all $R\ge1$, $0\le\gamma<\beta$, $0<\alpha-d\ll1$, and $0<p-1\ll_{\alpha}1$,
\begingroup\allowdisplaybreaks
\begin{eqnarray}
\lefteqn{\quad\Norm{\partial^\fun E(f,g)}_\beta^2\,\lesssim_{\alpha,p}\,R^d\pi_*(R)^{-1}\mathrm{[RHS\eqref{est2-LSI}]}}\nonumber\\
&&+\,R^{-\beta}\Big( \int \calC^2\mu_*^2[\nabla g]_\infty^2\Big)\bigg(\Big( \int |f|^2\Big)+(Rr_*(0))^{d\frac{p-1}{p}}\Big( \int [\nabla\bar u]_\infty^{2p}\Big)^\frac{1}p\bigg)\nonumber\\
&&+\,R^{-\beta}\Big( \int \calC^2\mu_*^2[\nabla^2\bar u]_\infty^2\Big)\bigg(\Big( \int r_*^d[g]_\infty^2\Big)+(Rr_*(0))^{d\frac{p-1}{p}}\Big( \int [g]_\infty^{2p}\Big)^\frac{1}p\bigg)\nonumber\\
&&+\,R^{d-\beta}r_*(0)^{d-\gamma+\alpha\frac{p-1}{2p}}\Big( \int r_*^{\gamma\frac{2p}{p-1}}w_R^{-\alpha}\Big)^\frac{p-1}{2p}\nonumber\\
&&\qquad\times\bigg(\Big( \int \calC^4\mu_*^4[\nabla g]_\infty^4\Big)^\frac{1}{2}\Big( \int w_R^{2p(d-\gamma)+\alpha(p-1)}(|f|+[\nabla\bar u]_\infty)^{4p}\Big)^\frac{1}{2p}\nonumber\\
&&\qquad\qquad\quad+\,\Big( \int \calC^4\mu_*^4[\nabla^2\bar u]_\infty^4\Big)^\frac{1}{2}\Big( \int r_*^{2pd}w_R^{2p(d-\gamma)+\,\alpha(p-1)}[g]_\infty^{4p}\Big)^\frac{1}{2p}\bigg),\label{est2-WLSI}
\end{eqnarray}
\endgroup
where we use the short-hand notation $\mathrm{[RHS\eqref{est2-LSI}]}$ for the RHS of~\eqref{est2-LSI}.\qedhere
\end{enumerate}
\end{prop}
The proofs of Propositions~\ref{prop:boundG2} and~\ref{prop:boundG1} rely on two further ingredients: large-scale weighted Calder\'on-Zygmund estimates and moment bounds on the extended corrector $(\phi,\sigma)$ (which are at the origin of the scaling $\mu_{*}$ in the estimates).
We start by recalling the former, which follows from~\cite[Theorem~1, Corollary~4, and Corollary~5]{GNO-reg} (see also~\cite[Section~7]{AKM-book}).
\begin{lem}[\cite{GNO-reg}]\label{lem:r*}
Assume that the coefficient field $\Aa$ is Gaussian with parameter $\beta>0$, and let $\pi_*$ be as in~\eqref{e.pi*}.
There exists a stationary, $\frac 18$-Lipschitz continuous random field $r_*\ge 1$ (the so-called \emph{minimal radius}), satisfying for some (deterministic) constant $C\simeq1$,
\begin{equation}\label{e.int-sto-r*}
\expec{\exp\Big(\frac1C \pi_*(r_*)\Big)}\le2,
\end{equation}
such that the following properties hold a.s.,
\begin{enumerate}[(a)]
\item \emph{Mean-value property:}\\For any $\Aa$-harmonic function $u$ in $B_R$ (that is, $-\nabla \cdot \Aa\nabla u=0$ in $B_R$), we have for all radii $r_*(0)\leq r\leq R$,
\begin{align}\label{eq:mean-value}
\fint_{B_r} |\nabla u|^2 \,\lesssim\, \fint_{B_R}|\nabla u|^2.
\end{align}
Applied to the extended corrector of Lemma~\ref{si}, this yields for all $\ell \ge 1$
and $x\in \R^d$,
\begin{align}\label{eq:mean-value-corr}
\int_{B_\ell(x)} |\nabla (\phi,\sigma)|^2 \,\lesssim\, \big(\ell+r_*(x)\big)^d.
\end{align}
\item 
 \emph{Large-scale Calder\'on-Zygmund estimates:}\\Set $B_*(x):=B_{r_*(x)}(x)$, and
 more generally $B_{\ell*}(x):=B_{\ell+r_*(x)}(x)$.
For all $1<p<\infty$, for all (sufficiently fast) decaying scalar fields $u$ and vector fields $g$ related in $\R^d$ by
\begin{equation*} 
-\nabla \cdot a \nabla u=\nabla \cdot g,
\end{equation*}
we have
\begin{equation}\label{I1-no-weight}
\int\Big(\fint_{B_*(x)}|\nabla u|^2\Big)^\frac{p}{2} dx\,\lesssim_p\, \int\Big(\fint_{B_*(x)}|g|^2\Big)^\frac{p}{2} dx.
\end{equation}
\item[(c)] \emph{Large-scale weighted Calder\'on-Zygmund estimates:}\\For all $2\le p <\infty$, $0\le \gamma < d(p-1)$, and for all non-decreasing radial weights $w\ge1$ satisfying
\begin{equation*}
w(r)\le w(r')\le \Big(\frac{r'}r\Big)^\gamma w(r)\qquad\text{for all $0\le r\le r'$},
\end{equation*}
we have for all $u$ and $g$ as in~(b) above,
\begin{equation}\label{I1-weight}
\int\Big(\fint_{B_*(x)}|\nabla u|^2\Big)^\frac{p}{2}w_*(x)\,dx
\,\lesssim_{p,\gamma}\, 
 \int\Big(\fint_{B_*(x)}|g|^2\Big)^\frac{p}{2} w_*(x)\,dx, 
\end{equation}
where $w_*(x):=w(|x|+r_*(0))$.
\qedhere
\end{enumerate}
\end{lem}
Whereas the minimal radius $r_*$ quantifies the sublinearity of the extended corrector at infinity~\cite{GNO-reg}, the precise growth of the latter is estimated as follows (cf.~\cite[Theorem~2]{GNO-quant}).
\begin{lem}[\cite{GNO-reg,GNO-quant}]\label{lem:bd-sigma}
Assume that the coefficient field $\Aa$ is Gaussian with parameter $\beta>0$, let $\mu_*$ be as
in~\eqref{mudbeta}, and let $r_*$ be as in Lemma~\ref{lem:r*}.
Then the extended corrector $(\phi,\sigma)$ defined in Lemma~\ref{si} satisfies for all $x\in \R^d$,
\begin{equation}\label{eq:pr-corr-4}
\Big(\fint_{B(x)} |(\phi,\sigma)|^2\Big)^{\frac{1}{2}}\,\le\,\calC(x) \mu_{*}(x),
\end{equation}
where $\calC\ge1$ is a 1-Lipschitz continuous random field with stretched exponential moments: there exist $\gamma\simeq_{\beta}1$ and $C_\gamma\simeq_\gamma1$
such that
\begin{equation}\label{e.int-sto-C}
\expec{\exp\Big(\frac1{C_\gamma}\calC^{\gamma}\Big)}\,\le\,2.\qedhere
\end{equation}
\end{lem}
In order to reformulate integrals in a form well-suited to apply (weighted) large-scale Calder\'on-Zygmund estimates,
we display below an auxiliary lemma that takes advantage of the Lipschitz continuity of $r_*$.
\begin{lem}\label{lem:integrals}
Let $\|\cdot\|_\ell$ be defined in~\eqref{eq:partition:0} and let $r_*$ be as in Lemma~\ref{lem:r*}.
For all $U,V$ and $\ell\ge 1$, we have
\begin{equation}\label{e.int-r_*0}
\|UV\|_\ell^2 \,\lesssim\,\int\Big(\int_{B_{2\ell*}(x)} |U|^2\Big) \Big(\fint_{B_*(x)} |V|^2\Big)\,dx,
\end{equation}
and the refined estimate
\begin{multline}\label{e.int-r_*}
\|UV\|_\ell^2 \,\lesssim\,\int_{|x|\ge \ell} \Big(\int_{B_{2\ell*}(x)} |U|^2\Big) \Big(\fint_{B_*(x)} |V|^2\Big)\,dx\\
+\bigg(\int_{B_{7\ell*}(0)} \Big(\fint_{\bar B_*(x)}|U|^2\Big)^\frac12 \Big(\fint_{\bar B_*(x)}|V|^2\Big)^\frac12\,dx\bigg)^2,
\end{multline}
where we recall $B_*(x)=B_{r_*(x)}(x)$, $B_{2\ell*}(x)=B_{2\ell+r_*(x)}(x)$, and $B_{7\ell*}(0)=B_{7\ell+r_*(0)}(0)$, and where we have set $\bar B_*(x):=B_{5r_*(x)} (x)$.
\end{lem}

\section{Proof of the representation formulas and of the main estimates} \label{sec:pr-main}

\subsection{Proof of Lemma~\ref{lem:decompI3eps}: Representation formulas}

We first introduce some notation.
Let $\Aa$ and $\tilde\Aa$ be two (admissible) coefficient fields, and set $\delta \Aa:=\tilde\Aa-\Aa$.
For all random variables (or fields) $F=F(\Aa)$, we set $\tilde F:=F(\tilde\Aa)$ and $\delta F:=\tilde F-F$.
We then denote by $(\phi,\sigma)$, $(\phi^*,\sigma^*)$, $(\tilde\phi,\tilde\sigma)$, and $(\tilde\phi^*,\tilde\sigma^*)$ the extended correctors associated with $\Aa$, $\Aa^*$, $\tilde\Aa$, and $\tilde\Aa^*$, respectively.

\medskip

\step1 Proof of identity~\eqref{eq:decomp-der-J0}.
\nopagebreak

\noindent The definition \eqref{intro:HC} of $\Xi_{ij}$ yields
\[\delta J_0(\F)= \int \F_{ij}\delta \Xi_{ij}= \int \F_{ij}\,e_j\cdot (\Aa-\bar\Aa)\nabla \delta \phi_i + \int \F_{ij}\,e_j \cdot \delta \Aa (\nabla \tilde\phi_i+e_i).\]
Using the definition~\eqref{eq:aux-seps} of the auxiliary field $S$ as well as the corrector equation~\eqref{f.2} for $\phi_i$ and $\tilde\phi_i$ in the form
\begin{align}\label{eq:corr-eqn-der}
\nabla\cdot\Aa\nabla\delta\phi_i=-\nabla\cdot\delta\Aa(\nabla\tilde\phi_i+e_i),
\end{align}
the first RHS term above can be rewritten as
\[ \int \F_{ij}\,e_j\cdot (\Aa-\bar\Aa)\nabla \delta \phi_i=- \int \nabla S_i\cdot\Aa\nabla\delta \phi_i= \int \nabla S_i\cdot\delta\Aa(\nabla\tilde\phi_i+e_i),\]
and the conclusion~\eqref{eq:decomp-der-J0} follows from the definition~\eqref{s.50} of the functional derivative.

\medskip
\step2 Proof of identity~\eqref{eq:decomp-der-E}.
\nopagebreak

\noindent We start by giving a suitable representation formula for the functional derivative of the homogenization commutator of the solution $(\Aa-\bar\Aa)\nabla u$.
By the property~\eqref{f.5} of the flux corrector $\sigma_j^*$ in the form $(\Aa^*-\bar\Aa^*)e_j=-\Aa^*\nabla\phi_j^*+\nabla\cdot\sigma_j^*$, and by its skew-symmetry~\eqref{f.19} in the form $(\nabla\cdot\sigma_j^*)\cdot\nabla \delta u=-\nabla\cdot(\sigma_j^*\nabla \delta u)$, we find
\begin{eqnarray*}
\delta\big(\ee_j\cdot(\Aa-\bar\Aa)\nabla u\big)&=&\ee_j\cdot\delta\Aa\nabla\tilde u+\ee_j\cdot(\Aa-\bar\Aa)\nabla\delta u\\
&=&\ee_j\cdot\delta\Aa\nabla\tilde u-\nabla\cdot(\sigma_j^*\nabla \delta u)-\nabla\phi_j^*\cdot\Aa\nabla\delta u.
\end{eqnarray*}
Equation~\eqref{eq:first-def-ups} for $u$ and $\tilde u$ in the form
\begin{align}\label{eq:u-eqn-der}
-\nabla\cdot\Aa\nabla\delta u=\nabla\cdot\delta\Aa\nabla\tilde u
\end{align}
allows us to rewrite the last RHS term as
\begin{eqnarray*}
-\nabla\phi_j^*\cdot\Aa\nabla\delta u&=&-\nabla\cdot(\phi_j^*\Aa\nabla\delta u)+\phi_j^*\nabla\cdot\Aa\nabla\delta u\\
&=&-\nabla\cdot(\phi_j^*\Aa\nabla\delta u)-\phi_j^*\nabla\cdot\delta\Aa\nabla\tilde u\\
&=&-\nabla\cdot(\phi_j^*\Aa\nabla\delta u)-\nabla\cdot(\phi_j^*\delta\Aa\nabla\tilde u)+\nabla\phi_j^*\cdot\delta\Aa\nabla\tilde u.
\end{eqnarray*}
Hence, we conclude
\[\delta\big(\ee_j\cdot(\Aa-\bar\Aa)\nabla u\big) = (\nabla\phi_j^*+\ee_j)\cdot\delta\Aa\nabla\tilde u-\nabla\cdot\big((\phi_j^*\Aa+\sigma_j^*)\nabla \delta u\big)-\nabla\cdot(\phi_j^*\delta\Aa\nabla\tilde u),\]
and similarly, replacing $x\mapsto u(x)$ by $x\mapsto\phi_i(x)+x_i$,
\[\delta\Xi_{ij} = (\nabla\phi_j^*+\ee_j)\cdot\delta\Aa(\nabla\tilde\phi_i+\ee_i)-\nabla\cdot\big((\phi_j^*\Aa+\sigma_j^*)\nabla \delta\phi_i\big)-\nabla\cdot\big(\phi_j^*\delta\Aa(\nabla\tilde\phi_i+\ee_i)\big).\]
Considering $\delta\big(\ee_j\cdot(\Aa-\bar\Aa)\nabla u\big)-\nabla_i\bar u\,\delta\Xi_{ij}$ and multiplying by $g_j$,
we are led to
\begin{multline*}
\delta E(f,g)= \int g_j(\nabla\phi_j^*+\ee_j)\cdot\delta\Aa\big(\nabla\tilde u-(\nabla\tilde\phi_i+\ee_i)\nabla_i\bar u\big)\\
+ \int \phi_j^*\nabla g_j\cdot\delta\Aa\nabla\tilde u- \int \phi_j^*\nabla(g_j\nabla_i\bar u)\cdot\delta\Aa(\nabla\tilde\phi_i+\ee_i)\\
+ \int \nabla g_j\cdot(\phi_j^*\Aa+\sigma_j^*)\nabla \delta u- \int \nabla(g_j\nabla_i\bar u)\cdot(\phi_j^*\Aa+\sigma_j^*)\nabla \delta\phi_i.
\end{multline*}
For the first RHS term we use the definition of $w_f$ in form of $\nabla u-(\nabla\phi_i+\ee_i)\nabla_i\bar u=\nabla w_f+\phi_i\nabla\nabla_i\bar u$, 
whereas for the last two RHS terms we use the definitions~\eqref{eq:aux-reps} and~\eqref{eq:aux-Reps} of the auxiliary fields $r$ and $R$, combined with equations~\eqref{eq:corr-eqn-der} and~\eqref{eq:u-eqn-der}, so that
\begin{multline}\label{eq:corr-eqn-der-pre}
\delta E(f,g)= \int g_j(\nabla\phi_j^*+\ee_j)\cdot\delta\Aa(\nabla\tilde w_f+\tilde\phi_i\nabla\nabla_i\bar u)\\
+ \int \phi_j^*\nabla g_j\cdot\delta\Aa\nabla\tilde u- \int \phi_j^*\nabla(g_j\nabla_i\bar u)\cdot\delta\Aa(\nabla\tilde\phi_i+\ee_i)\\
+ \int \nabla r\cdot\delta\Aa\nabla\tilde u- \int \nabla R_{i}\cdot\delta\Aa(\nabla\tilde\phi_i+\ee_i),
\end{multline}
and the conclusion~\eqref{eq:decomp-der-E} follows from the definition~\eqref{s.50} of the functional derivative.


\subsection{Proof of Lemma~\ref{lem:integrals}}\label{sec:tria}

We first recall the following equivalence for all non-negative functions $h$,
\begin{equation}\label{e.equiv-B*}
\int h\simeq\int \Big(\fint_{B_*(x)}h\Big),
\end{equation}
cf.~\cite[Proof of Corollary~4, Step~5]{GNO-reg}.
Estimate~\eqref{e.int-r_*0} is a consequence of \eqref{e.equiv-B*} in form of 
\begin{eqnarray*}
 \|UV\|_\ell^2\,=\,\ell^{-d}\int\Big(\int_{B_\ell(x)} |UV|\Big)^2 dx&\le& \int \Big(\int_{B_\ell(x)}|U|^2\Big)\Big(\fint_{B_\ell(x)}|V|^2\Big)\,dx
\nonumber
\\
&\le &\int \bigg(\fint_{B_\ell(x)} |V(y)|^2\Big(\int_{B_{2\ell}(y)} |U|^2\Big)dy\bigg)\,dx\nonumber \\
&=&\int |V(x)|^2\Big(\int_{B_{2\ell}(x)} |U|^2\Big)\,dx \nonumber \\
& \lesssim&\int\Big(\int_{B_{2\ell*}(x)} |U|^2\Big) \Big(\fint_{B_*(x)} |V|^2 \Big)\,dx.
\end{eqnarray*}
We now turn to the proof of~\eqref{e.int-r_*}.
We distinguish the generic case $r_*(0)\le \ell$ from the non-generic case $r_*(0)>\ell$, and we start with the latter.
By the $\frac18$-Lipschitz continuity of~$r_*$ and the assumption $r_*(0)>\ell\ge 1$, we have  
$$|x|\le \ell\quad \implies \quad \tfrac78 r_*(0) \le r_*(0)-\tfrac18\ell \le r_*(x)\le r_*(0)+\tfrac 18\ell \le \tfrac98 r_*(0),$$
and similarly,
$$|x|\le \tfrac \ell 2 \quad \implies \quad \tfrac{15}{16} r_*(0)\le r_*(x) \le \tfrac{17}{16}r_*(0) \quad \implies \quad B_{4r_*(0)}(0)\subset \bar B_*(x),$$
since $\frac{15}{16}\times 5-\frac12=\frac{67}{16}>4$, where we recall the definition $\bar B_*(x)=B_{5r_*(x)}(x)$. Hence,
\begin{eqnarray*}
\lefteqn{\int_{|x|\le \ell}\Big(\int_{B_{2\ell*}(x)} |U|^2\Big) \Big(\fint_{B_*(x)} |V|^2 \Big)\,dx}\\
&\lesssim&r_*(0)^{-d}\,\ell^d  \Big(\int_{B_{(3+\frac98)r_*(0)}(0)} |U|^2\Big)\Big(\int_{B_{(1+\frac98)r_*(0)}(0)} |V|^2\Big) 
\\
&\le& r_*(0)^{-d}\,   \bigg(\int_{B_{\frac \ell2}(0)} \Big(\int_{\bar B_{*}(x)} |U|^2\Big)^\frac12\Big(\int_{\bar B_{*}(x)} |V|^2\Big)^\frac12\bigg)^2 
\\
&\lesssim& \bigg(\int_{B_{\ell*}(0)} \Big(\fint_{\bar B_{*}(x)} |U|^2\Big)^\frac12\Big(\fint_{\bar B_{*}(x)} |V|^2\Big)^\frac12\bigg)^2.
\end{eqnarray*}
Combined with~\eqref{e.int-r_*0}, this yields the conclusion~\eqref{e.int-r_*} for $r_*(0)>\ell$.
We turn to the generic case $r_*(0)\le \ell$.
We split the integral over $\R^d$ into the far-field contribution $|x|\ge 4\ell$ and the near-field contribution $|x|<4\ell$. 
For the former, we proceed as above,
\begin{eqnarray*}
\ell^{-d} \int_{|x|\ge 4\ell}    \Big(\int_{B_\ell(x)}|UV|\Big)^2dx
&\le&\int_{|x|\ge 4\ell} \bigg(\fint_{B_\ell(x)} |V(y)|^2  \Big(\int_{B_{2\ell}(y)} |U|^2\Big)\,dy  \bigg)\,dx  \\
&\le&\int_{|x|\ge3 \ell}  |V(x)|^2 \Big(\int_{B_{2\ell}(x)} |U|^2\Big)\,dx\\
&\lesssim&\int \bigg(\fint_{B_*(x)}|V(y)|^2\Big(\int_{B_{2\ell}(y)} |U|^2\Big) \mathds 1_{|y|\ge 3 \ell}\,dy\bigg) dx,
\end{eqnarray*}
where the last bound follows from~\eqref{e.equiv-B*}.
By the $\frac18$-Lipschitz continuity of $r_*$ and the assumption $r_*(0)\le \ell$,
we infer that the condition $|x|<\ell$ implies $r_*(x) \le r_*(0)+\frac18\ell < 2\ell$, hence $B_*(x) \subset B_{3\ell}(0)$. The above inequality then reduces to 
\begin{eqnarray}
\ell^{-d} \int_{|x|\ge 4\ell}\Big(\int_{B_\ell(x)}|UV|\Big)^2dx
&\lesssim& \int_{|x|\ge \ell}\bigg(\fint_{B_*(x)}|V(y)|^2\Big(\int_{B_{2\ell}(y)} |U|^2\Big)\,dy\bigg)\,dx\nonumber\\
&\le & \int_{|x|\ge \ell} \Big(\int_{B_{2\ell*}(x)} |U|^2\Big) \Big(\fint_{B_*(x)}|V|^2\Big)\,dx.\label{e.a.+1}
\end{eqnarray}
We turn to the near-field contribution $|x|<4\ell$. We start with the trivial estimate
$$
\ell^{-d} \int_{|x|< 4\ell}\Big(\int_{B_\ell(x)} |UV|\Big)^2dx \,\lesssim \,\Big(\int_{B_{5\ell}(0)} |UV| \Big)^2,
$$
and we use~\eqref{e.equiv-B*} in form of
$$
\int_{B_{5\ell}(0)} |UV| \lesssim \int \Big(\fint_{B_*(x)} |UV|\mathds 1_{|y|<5\ell} \,dy\Big)\,dx.
$$
By the $\frac18$-Lipschitz continuity of $r_*$ and the assumption $r_*(0)\le \ell$, we infer that the condition $|x|>{7}\ell$ implies $B_*(x) \cap B_{5\ell}(0)=\varnothing$,
hence 
$$
\int_{B_{5\ell}(0)} |UV| \lesssim \int_{|x|<7\ell}\Big( \fint_{B_*(x)} |UV|\Big)\,dx.
$$
The Cauchy-Schwarz' inequality then leads to
$$
\ell^{-d} \int_{|x|< 4\ell}\Big(\int_{B_\ell(x)} |UV| \Big)^2dx \,\lesssim \,\bigg(\int_{|x|<7\ell} \Big(\fint_{B_*(x)} |U|^2\Big)^\frac12 \Big(\fint_{B_*(x)} |V|^2\Big)^\frac12dx\bigg)^2,
$$
and~\eqref{e.int-r_*} follows in combination with \eqref{e.a.+1} in the generic case $r_*(0)\le \ell$.


\subsection{Proof of Proposition~\ref{prop:boundG2}: Main estimates}\label{chap:boundG2}

We split the proof into two main steps, first addressing the case of the standard LSI ($\beta>d$), and then turning to the general multiscale case ($\beta\le d$). Let $R\ge 1$ be arbitrary.

\medskip

\step1
Proof of~\eqref{est1-LSI} for standard LSI ($\beta>d$).

\noindent Since for standard LSI ($\beta>d$) we have $\Norm{\partial^\fun J_0(\F)}_\beta\lesssim\|\partial^\fun J_0(\F)\|_{1}$ (cf.~\eqref{e.operator-norm-simpl}), it suffices to prove the following estimate: for all $\ell\ge 1$, $\alpha>d$, and $p>1$ with $\alpha(p-1)<d(2p-1)$,
\begin{equation}\label{est-G-LSI}
\|\partial^\fun J_0(\F)\|_{\ell}^2\,
\lesssim_{\alpha,p}\, \ell^dr_*(0)^{\alpha\frac{p-1}p}
\Big(\int r_*^{d\frac p{p-1}}w_{R}^{-\alpha}\Big)^\frac{p-1}p\Big(\int w_R^{\alpha(p-1)}|F|^{2p}\Big)^\frac1p.
\end{equation}
Starting point is formula~\eqref{eq:decomp-der-J0}, which, by~\eqref{e.int-r_*0} in Lemma~\ref{lem:integrals} 
for $U=\nabla \phi+\Id$ and $V=|F|+|\nabla S|$, and by the mean-value property~\eqref{eq:mean-value-corr}, implies
\begin{eqnarray}
\|\partial^\fun J_0(\F)\|_{\ell}^2 & \lesssim&\int \big(\ell+r_*(x)\big)^d\Big(\fint_{B_*(x)} |F|^2+|\nabla S|^2\Big)\, dx \nonumber
\\
&\lesssim &\ell^d \int r_*(x)^d \Big(\fint_{B_*(x)} |F|^2+|\nabla S|^2\Big)\, dx
.\label{e.decompJ_0-genell}
\end{eqnarray}
This yields the conclusion~\eqref{est-G-LSI} in combination 
with the following estimate applied for $s=q=1$ and $v=S$ (cf.~\eqref{eq:aux-seps}):
 If $v$ is the Lax-Milgram solution of $-\nabla\cdot\Aa\nabla v=\nabla\cdot h$ with $h\in C^\infty_c(\R^d)$, then for all $s\ge0$, $q\ge1$, $\alpha>d$, and $p>1$ with $\alpha(p-1)<d(2pq-1)$, 
\begin{multline}\label{eq-00-pre}
\int r_*(x)^{ds} \Big(\fint_{B_{*}(x)}|h|^2+|\nabla v|^2\Big)^qdx\\
\lesssim_{\alpha,p,q,s}\, r_*(0)^{\alpha\frac{p-1}p}\Big( \int r_*^{ds\frac p{p-1}}w_{R}^{-\alpha}\Big)^\frac{p-1}p\Big( \int w_R^{\alpha(p-1)}[h]_2^{2pq}\Big)^\frac1p,
\end{multline}
where for $s=0$ we may even choose $p=1$, in which case \eqref{eq-00-pre} is replaced by
\begin{align}\label{eq-00-pre-1}
\int\Big(\fint_{B_{*}(x)}|h|^2+|\nabla v|^2\Big)^qdx
\lesssim_q\,\int [h]_2^{2q},
\end{align}
which we state and prove here for future reference only.

\medskip\noindent
Here comes the argument for \eqref{eq-00-pre}.
For all $\alpha>d$ and $p>1$, we smuggle in the weight $w_{R*}(x):=w_R(|x|+r_*(0))$ to the power $\alpha\frac{p-1}{p}$, and use H\"older's inequality with exponents $(\frac{p}{p-1},p)$,
so that
\begin{multline*}
\int r_*(x)^{ds}\Big(\fint_{B_{*}(x)}|h|^2+|\nabla v|^2\Big)^qdx\,\lesssim_s\,
\Big(\int r_*^{ds\frac p{p-1}} w_{R*}^{-\alpha}\Big)^\frac{p-1}p\\
\times\bigg(\int w_{R*}(x)^{\alpha(p-1)}\Big(\fint_{B_*(x)}|h|^2+|\nabla v|^2\Big)^{pq}dx\bigg)^\frac1p,
\end{multline*}
where the first RHS sum is bounded by $ \int r_*^{ds\frac p{p-1}}w_{R}^{-\alpha}$ since $w_{R}\le w_{R*}$.
Provided that $\alpha(p-1)<d(2pq-1)$, we may apply the large-scale weighted Calder\'on-Zygmund estimate of~Lemma~\ref{lem:r*}(c) to the equation for $v$, to the effect of
\begin{multline*}
\int r_*(x)^{ds}\Big(\fint_{B_{*}(x)}|h|^2+|\nabla v|^2\Big)^qdx\,\lesssim_{\alpha,p,q,s}\,\Big(\int r_*^{ds\frac p{p-1}}w_{R}^{-\alpha}\Big)^\frac{p-1}p\\
\times\bigg(\int w_{R*}(x)^{\alpha(p-1)}\Big(\fint_{B_*(x)} |h|^2\Big)^{pq}dx\bigg)^\frac1p.
\end{multline*}
The claim~\eqref{eq-00-pre} then follows from the bound $w_{R*}(x)\le r_*(0)\inf_{B_*(x)}w_R$, Jensen's inequality, and~\eqref{e.equiv-B*}.
For $s=0$, we appeal to the large-scale (not weighted) Calder\'on-Zygmund estimate of~Lemma~\ref{lem:r*}(b), which amounts to choosing $p=1$ in the above.


\medskip

\step2 Proof of~\eqref{est1-WLSI} in the general multiscale case ($\beta\le d$).
\nopagebreak

\noindent The combination of~\eqref{e.operator-norm} with~\eqref{est-G-LSI} is not enough to prove~\eqref{est1-WLSI}, and we have 
to refine~\eqref{est-G-LSI} in the regime $\ell\ge R$.
By~\eqref{e.int-r_*} in Lemma~\ref{lem:integrals} and the mean-value property~\eqref{eq:mean-value-corr},
\begin{multline}\label{e.specific-L1-Linfty}
\|\partial^\fun J_0(\F)\|_{\ell}^2\,\lesssim\,\int_{|x|\ge \ell} \big(\ell+r_* (x)\big)^d \Big(\fint_{B_*(x)}{|F|^2+|\nabla S|^2}\Big)\, dx
\\
+ \bigg(\int_{B_{7\ell*}(0)} \Big(\fint_{\bar B_*(x)}|F|^2+|\nabla S|^2\Big)^\frac12\bigg)^2.
\end{multline}
Let now $v$ be the Lax-Milgram solution of $-\nabla\cdot\Aa\nabla v=\nabla\cdot h$ with $h\in C^\infty_c(\R^d)$.
In the following two substeps we estimate the far-field and near-field contributions separately.

\medskip

\substep{2.1} Far-field estimate: For all $s\ge0$, $q\ge1$, $0\le\gamma\le ds$, $\alpha>d$, and $p>1$ with $(ds-\gamma)p+\alpha(p-1)<d(2pq-1)$,
\begin{multline}\label{eq-00-pre+}
\int_{|x|\ge \ell} \big(\ell+r_*(x)\big )^{ds}\Big(\fint_{B_*(x)} |h|^2+|\nabla v|^2\Big)^qdx
\,\lesssim_{\gamma,\alpha,p,q,s}\, \ell^{\gamma}R^{ds-\gamma}r_*(0)^{ds-\gamma+\alpha\frac{p-1}p}\\
\times\Big( \int r_*^{\gamma\frac{p}{p-1}}w_R^{-\alpha}\Big)^\frac{p-1}p
\Big( \int w_R^{(ds-\gamma)p+ \alpha(p-1)}[h]_2^{2pq}\Big)^\frac1p.
\end{multline}
We smuggle in the weight $w_{R*}$ to the power $\alpha\frac{p-1}p$ and the weight $w_{1*}$ to the power $ds-\gamma$, and use H\"older's inequality with exponents~$(\frac{p}{p-1},p)$, to the effect of
\begin{multline*}
\int_{|x|\ge \ell} \big(\ell+r_*(x)\big)^{ds}
\Big(\fint_{B_*(x)}{|h|^2+|\nabla v|^2}\Big)^q dx\\
\lesssim_s\,
\bigg( \int_{|x|\ge \ell} w_{R*}^{-\alpha}w_{1*}^{-(ds-\gamma)\frac{p}{p-1}}(\ell+r_*)^{ds\frac{p}{p-1}}\bigg)^\frac{p-1}p
\\
\times\bigg(\int w_{R*}^{\alpha(p-1)} w_{1*}^{(ds-\gamma)p}\Big(\fint_{B_*(x)}|h|^2+|\nabla v|^2\Big)^{pq}dx\bigg)^\frac1p.
\end{multline*}
In the first RHS factor, we use the bound $w_{1*}(x) \gtrsim\ell+r_*(x)$ for $|x|\ge\ell$, while in the second RHS factor we use $w_{1*}\le Rw_{R*}$. The above then leads to
\begin{multline*}
\int_{|x|\ge \ell} \big(\ell+r_* (x)\big)^{ds}
\Big(\fint_{B_*(x)}{|h|^2+|\nabla v|^2}\Big)^q dx
\,\lesssim_{\gamma,s}\,R^{ds-\gamma}
\bigg( \int_{|x|\ge \ell} w_{R*}^{-\alpha}(\ell+r_*)^{\gamma\frac{p}{p-1}}\bigg)^\frac{p-1}p\\
\times\bigg(\int w_{R*}^{(ds-\gamma)p+ \alpha(p-1)}\Big(\fint_{B_*(x)}|h|^2+|\nabla v|^2\Big)^{pq}dx\bigg)^\frac1p.
\end{multline*}
Provided $(ds-\gamma)p+\alpha(p-1)<d(2pq-1)$, we may apply the large-scale weighted Calder\'on-Zygmund estimate of Lemma~\ref{lem:r*}(c) to the equation for $v$.
Using the bound $w_{R*}(x)\le r_*(0)\inf_{B_*(x)}w_R$, Jensen's inequality, and~\eqref{e.equiv-B*}, the conclusion~\eqref{eq-00-pre+} follows.

\medskip

\substep{2.2} Near-field estimate: For all $\ell\ge1$ and $p>1$,
\begin{eqnarray}
\int_{B_{7\ell*}(0)} \Big(\fint_{\bar B_*(x)}{|h|^2+|\nabla v|^2}\Big)^\frac12 
&\lesssim_p& \ell^{d\frac{p-1}p}r_*(0)^{d\frac{p-1}p} \ \Big( \int r_*^{d(1-\frac p2)_+}[h]_2^p\Big)^\frac1p.
\label{p2.2-2.1-term1}
\end{eqnarray}
Indeed, by H\"older's inequality with exponents~$(\frac p{p-1},p)$,
\begin{eqnarray*}
\lefteqn{\int_{B_{7\ell*}(0)} \Big(\fint_{\bar B_*(x)}|h|^2+|\nabla v|^2\Big)^\frac12}
\\
 &\lesssim & \ell^{d\frac{p-1}p}r_*(0)^{d\frac{p-1}p} \bigg(\int_{B_{7\ell*}(0)} \Big(\fint_{\bar B_*(x)}{|h|^2+|\nabla v|^2}\Big)^\frac p2\bigg)^\frac 1p
 \\
 &\lesssim& \ell^{d\frac{p-1}p}r_*(0)^{d\frac{p-1}p} \bigg(\int\Big(\fint_{B_*(x)}{|h|^2+|\nabla v|^2}\Big)^\frac p2\bigg)^\frac 1p,
\end{eqnarray*}
where we pass from $\bar B_*(x)=B_{5r_*(x)}(x)$ to $B_*(x)=B_{r_*(x)}(x)$.
The claim~\eqref{p2.2-2.1-term1} follows from the large-scale Calder\'on-Zygmund estimate of Lemma~\ref{lem:r*}(b) applied to $v$.

\medskip

\substep{2.3} Conclusion.

\noindent For all $\ell\ge1$, $0\le \gamma\le d$, $\alpha>d$, and $p>1$ with $(d-\gamma)p+\alpha(p-1)<d(2p-1)$ and $p\le2$, the combination of Substeps~2.1--2.2 with~\eqref{e.specific-L1-Linfty}
yields the following improvement of~\eqref{est-G-LSI},
\begin{multline}\label{est-G-WLSI}
\|\partial^\fun J_0(\F)\|_{\ell}^2\,\lesssim_{\gamma,\alpha,p}\,\ell^{2d\frac{p-1}p}
\Big( \int r_*^{d(1-\frac p2)}[F]_2^p\Big)^\frac2p\\
+\ell^{\gamma}R^{d-\gamma} r_*(0)^{d-\gamma+\alpha\frac{p-1}p}\Big( \int r_*^{\gamma\frac{p}{p-1}}w_R^{-\alpha}\Big)^\frac{p-1}p\Big( \int w_R^{(d-\gamma)p+ \alpha(p-1)}|F|^{2p}\Big)^\frac1p.
\end{multline}
We appeal to~\eqref{e.operator-norm}, which we combine wtih~\eqref{est-G-LSI} for $\ell \le R$
and with~\eqref{est-G-WLSI} for $\ell > R$. 
Provided that $0\le\gamma<\beta$, we compute
\begin{eqnarray*}
\int_1^R \ell^{d-1-\beta} d\ell 
&\lesssim &
\left\{
\begin{array}{lll}
R^{d-\beta} &:& \beta<d
\\
\log R&:&\beta=d
\end{array}
\right\}\,\simeq\, R^{d} \pi_*(R)^{-1},
\\
\int_R^\infty \ell^{-1-\beta} d\ell &\lesssim & R^{-\beta},
\\
R^{d-\gamma} \int_R^\infty \ell^{-1-\beta+\gamma }d\ell &\lesssim&R^{d-\gamma -(\beta-\gamma)}\,=\,R^{d-\beta},
\end{eqnarray*}
and the conclusion~\eqref{est1-WLSI} follows.


\subsection{Proof of Proposition~\ref{prop:boundG1}: Main estimates (cont'd)}

By~\eqref{eq:decomp-der-E} in Lemma~\ref{lem:decompI3eps}, we have $\partial^\fun E(f,g)=G_1+G_2+G_3$ with
\begin{gather*}
G_1\,:=\,g_{j}\,(\nabla\phi_j^*+\ee_j)\otimes(\nabla w_{f}+\phi_i\nabla\nabla_i\bar u),\qquad G_2\,:=\,(\phi_j^*\,\nabla g_{j}+\nabla r)\otimes\nabla u,\\
G_3\,:=\,-\big(\phi_j^*\,\nabla(g_{j}\nabla_i\bar u)+\nabla R_{i}\big)\otimes(\nabla\phi_i+\ee_i),
\end{gather*}
so that it suffices to estimate the norms of each of the $G_i$'s separately.
We split the proof into two main steps: we first address the case of the standard LSI ($\beta>d$), and then turn to the general multiscale case ($\beta\le d$). Let $R\ge 1$ be arbitrary.

\medskip

\step1 Proof of~\eqref{est2-LSI} for standard LSI ($\beta>d$).\\
Since for standard LSI ($\beta>d$) we have $\Norm{\partial^\fun J_0(\F)}_\beta\lesssim\|\partial^\fun J_0(\F)\|_{1}$, it suffices to establish the following estimates: for all $\ell\ge 1$, $0<\alpha-d\ll1$, and $0<p-1\ll_\alpha 1$,
\begingroup\allowdisplaybreaks
\begin{eqnarray}
\|\G_1\|_{\ell}^2&\lesssim_{\alpha,p}&\ell^{d}\Big( \int r_*^{4d}[g]_\infty^4\Big)^\frac12\Big( \int \calC^4\mu_*^4[\nabla^2\bar u]_\infty^4\Big)^\frac12,\label{est1-G1}\\
\|\G_2\|_{\ell}^2&\lesssim_{\alpha,p}&\ell^{d}r_*(0)^{\alpha\frac{p-1}{2p}}\Big( \int r_*^{d\frac{2p}{p-1}}w_R^{-\alpha}\Big)^\frac{p-1}{2p}\nonumber\\
&&\qquad\qquad\times\Big( \int \calC^4\mu_*^4[\nabla g]_\infty^4\Big)^\frac12\Big( \int w_R^{\alpha(p-1)}|f|^{4p}\Big)^\frac1{2p},\label{est1-G2}\\
\|\G_3\|_{\ell}^2&\lesssim_{\alpha,p}&\ell^dr_*(0)^{\alpha\frac{p-1}p}\Big( \int r_*^{d\frac{p}{p-1}}w_R^{-\alpha}\Big)^\frac{p-1}p\nonumber\\
&&\qquad\qquad\times\Big( \int w_R^{\alpha(p-1)}\calC^{2p}\mu_*^{2p}[\nabla(g\nabla\bar u)]_\infty^{2p}\Big)^\frac1p.\label{est1-G3}
\end{eqnarray}
\endgroup
Indeed, replacing $p$ by $\frac{2p}{p+1}$, estimating $[\nabla(g\nabla\bar u)]_\infty\lesssim [\nabla g]_\infty[\nabla\bar u]_\infty+[g]_\infty[\nabla^2\bar u]_\infty$, and using H\"older's inequality, the estimate~\eqref{est1-G3} easily leads to
\begin{multline*}
\|\G_3\|_{\ell}^2\,\lesssim_{\alpha,p}\,\ell^dr_*(0)^{\alpha\frac{p-1}{2p}}\Big( \int r_*^{d\frac{2p}{p-1}}w_R^{-\alpha}\Big)^\frac{p-1}{2p}\\
\times\bigg(\Big( \int \calC^{4}\mu_*^4[\nabla^2\bar u]_\infty^4\Big)^\frac{1}{2}\Big( \int w_R^{\alpha(p-1)}[g]_\infty^{4p}\Big)^\frac{1}{2p}\\
+\Big( \int \calC^{4}\mu_*^4[\nabla g]_\infty^4\Big)^\frac{1}{2}\Big( \int w_R^{\alpha(p-1)}[\nabla\bar u]_\infty^{4p}\Big)^\frac{1}{2p}\bigg),
\end{multline*}
so that \eqref{est2-LSI} follows in combination with \eqref{est1-G1} and \eqref{est1-G2}.
We address the estimates~\eqref{est1-G1}--\eqref{est1-G3} separately, and split the proof into three substeps.

\medskip

\substep{1.1} Proof of~\eqref{est1-G1}.

\noindent 
By~\eqref{e.int-r_*0} in Lemma~\ref{lem:integrals},
\begin{equation*}
\|G_1\|_{\ell}^2\,\lesssim\,\int\Big(\int_{B_{2\ell *}(x)} |g|^2|\nabla\phi+\Id|^2\Big) \Big(\fint_{B_*(x)} (|\nabla w_f|+|\phi||\nabla^2\bar u|)^2 \Big)dx
\end{equation*}
which by Cauchy-Schwarz' inequality turns into
\begin{multline}\label{est1-G3-0}
\|G_1\|_{\ell}^2\,\lesssim\,
 \bigg(\int\Big(\int_{B_{2\ell*}(x)}|g|^2|\nabla\phi+\Id|^2\Big)^2dx\bigg)^\frac12\\
 \times\bigg(\int \Big(\fint_{B_*(x)}(|\nabla w_f|+|\phi||\nabla^2\bar u|)^2\Big)^2dx\bigg)^\frac12.
\end{multline}
We start by treating the first RHS factor. Taking the local supremum of $g$, using the 
mean-value property~\eqref{eq:mean-value-corr}, the Lipschitz continuity of $r_*$, Jensen's inequality, and~\eqref{e.equiv-B*}, we obtain
\begin{eqnarray}
{\int\Big(\int_{B_{2\ell*}(x)}|g|^2|\nabla\phi+\Id|^2\Big)^2dx}
&\lesssim& \int\Big(\int_{B_{2\ell*}(x)} r_*^d [g]_\infty^2 \Big)^2dx \nonumber
\\
&\lesssim&\int\Big(\fint_{B_{2\ell*}(x)}r_*^d(\ell+r_*)^d [g]_\infty^2 \Big)^2dx\nonumber
\\
&\le& \int r_*^{2d}(\ell+r_*)^{2d} [g]_\infty^4\nonumber
\\
&\lesssim&\ell^{2d} \int r_*^{4d}[g]_\infty^4.\label{est1-G3-01}
\end{eqnarray}
We turn to the second RHS factor in~\eqref{est1-G3-0}.
Note that the two-scale expansion error $w_f$ satisfies the following equation (cf.~\cite[proof of Theorem~3]{GNO-quant}),
\begin{align}\label{eq:wf}
-\nabla\cdot\Aa\nabla w_f=\nabla\cdot \big((\Aa\phi_j+\sigma_j)\nabla\nabla_j\bar u\big).
\end{align}
By~\eqref{eq-00-pre-1} with $q=2$ applied to $w_f$, we obtain after taking local suprema of $\nabla^2\bar u$, and controlling correctors by Lemma~\ref{lem:bd-sigma},
\begin{align}\label{est1-G3-02}
\int \Big(\fint_{B_*(x)}(|\nabla w_f|+|\phi||\nabla^2\bar u|)^2\Big)^2
\,\lesssim\, \int [(|\phi|+|\sigma|)\nabla^2\bar u]_2^4\,\lesssim\, \int \calC^4\mu_*^4[\nabla^2\bar u]_\infty^4.
\end{align}
Combined with~\eqref{est1-G3-0} and~\eqref{est1-G3-01}, this yields the conclusion~\eqref{est1-G1}.

\medskip

\substep{1.2} Proof of~\eqref{est1-G2}.

\noindent By Lemma~\ref{lem:integrals} in form of \eqref{e.int-r_*0}
and Cauchy-Schwarz' inequality,
\begin{align}\label{est1-G4-00}
\|G_2\|_{\ell}^2\,\lesssim\,\bigg(\int\Big(\fint_{B_*(x)}(|\phi||\nabla g|+|\nabla r|)^2\Big)^2\bigg)^\frac12\bigg(\int\Big(\int_{B_{2\ell*}(x)}|\nabla u|^2\Big)^2\bigg)^\frac12.
\end{align}
We start with the first RHS factor.
By~\eqref{eq-00-pre-1} with $q=2$ applied to the solution $r$ of~\eqref{eq:aux-reps}, we obtain after taking local suprema of $\nabla g$ and controlling
correctors by Lemma~\ref{lem:bd-sigma},
\begin{align}\label{est1-G4-01}
 \int\Big(\fint_{B_*(x)}(|\phi||\nabla g|+|\nabla r|)^2\Big)^2\,\lesssim\, \int [(|\phi|+|\sigma|)\nabla g]_2^4\,\lesssim\, \int \calC^4\mu_*^4[\nabla g]_\infty^4.
\end{align}
We turn to the second RHS factor in~\eqref{est1-G4-00}.
By the Lipschitz continuity of $r_*$, Jensen's inequality, and~\eqref{e.equiv-B*},
\begin{equation*}
\int\Big(\int_{B_{2\ell*}(x)}|\nabla u|^2\Big)^2\,\lesssim\, \int\Big(\fint_{B_{2\ell*}(x)}(\ell+r_*)^d |\nabla u|^2\Big)^2
\,\lesssim\,\int (\ell+r_*)^{2d} \Big(\fint_{B_{*}(x)} |\nabla u|^2\Big)^2.
\end{equation*}
By~\eqref{eq-00-pre} with $s=q=2$ applied to the solution $u$ of~\eqref{eq:first-def-ups}, we deduce for all $\alpha>d$ and $p>1$ with~\mbox{$\alpha(p-1)<d(4p-1)$},
\begin{multline}\label{eq:est-nabu-3}
\int (\ell+r_*)^{2d} \Big(\fint_{B_*(x)}|\nabla u|^2\Big)^2\\
\lesssim_{\alpha,p}\,\ell^{2d}r_*(0)^{\alpha\frac{p-1}p}\Big( \int r_*^{d\frac{2p}{p-1}}w_R^{-\alpha}\Big)^\frac{p-1}p\Big( \int w_R^{\alpha(p-1)}|f|^{4p}\Big)^\frac1p.
\end{multline}
Combined with~\eqref{est1-G4-00} and~\eqref{est1-G4-01}, this yields~\eqref{est1-G2}.

\medskip

\substep{1.3} Proof of~\eqref{est1-G3}.

\noindent
By~\eqref{e.int-r_*0} in Lemma~\ref{lem:integrals} and the mean-value property~\eqref{eq:mean-value-corr}, we find
\begin{eqnarray*}
\|G_3\|_{\ell}^2&\lesssim &
\int \big(\ell+r_*(x)\big)^d \Big(\fint_{B_*(x)}\big(|\phi^*||\nabla(g\nabla\bar u)|+|\nabla R|\big)^2\Big)\,dx.
\end{eqnarray*}
By~\eqref{eq-00-pre} with $s=q=1$ applied to the solution $R$ of~\eqref{eq:aux-Reps}, we deduce for all $\alpha>d$ and $p>1$ with $\alpha(p-1)<d(2p-1)$,
\begin{align*}
\|G_3\|_{\ell}^2\,\lesssim\,\ell^dr_*(0)^{\alpha\frac{p-1}p}\Big( \int r_*^{d\frac{p}{p-1}}w_R^{-\alpha}\Big)^\frac{p-1}p\Big( \int w_R^{\alpha(p-1)}[\phi^*\nabla(g\nabla\bar u)]_2^{2p}\Big)^\frac1p.
\end{align*}
Taking local suprema of $\nabla(g\nabla\bar u)$ and using Lemma~\ref{lem:bd-sigma} to control correctors,~\eqref{est1-G3} follows.

\medskip

\step2 Proof of~\eqref{est2-WLSI} in the general multiscale case ($\beta\le d$).

\noindent As in Step~2 of the proof of Proposition~\ref{prop:boundG2}, we need to refine~\eqref{est1-G1}--\eqref{est1-G3} in the range $\ell>R$. 
More precisely, we shall establish that for all $\ell\ge 1$, $0\le\gamma\le d$, $0<\alpha-d\ll1$, and $0<p-1\ll_\alpha1$,
\begingroup\allowdisplaybreaks
\begin{eqnarray}
\|G_1\|_{\ell}^2&\lesssim_{\gamma,\alpha,p}&\Big( \int \calC^2\mu_*^2[\nabla^2\bar u]_\infty^2\Big)\Big( \int r_*^d[g]_\infty^2\Big)\nonumber\\
\hspace{-0.7cm}&&+\ell^{\gamma}R^{d-\gamma}r_*(0)^{d-\gamma+\alpha\frac{p-1}{2p}}\Big( \int r_*^{\gamma\frac{2p}{p-1}}w_R^{-\alpha}\Big)^\frac{p-1}{2p}\nonumber\\
\hspace{-0.7cm}&&\hspace{1cm}\times\Big( \int r_*^{2pd}w_R^{2p(d-\gamma)+ \alpha(p-1)}[g]_\infty^{4p}\Big)^\frac1{2p}\Big( \int \calC^4\mu_*^4[\nabla^2\bar u]_\infty^4\Big)^\frac12,\label{est2-G1}\\
\|G_2\|_{\ell}^2&\lesssim_{\gamma,\alpha,p}&\Big( \int \calC^2\mu_*^2[\nabla g]_\infty^2\Big)\Big( \int |f|^2\Big)\nonumber\\
\hspace{-0.7cm}&&+\ell^{\gamma}R^{d-\gamma}r_*(0)^{d-\gamma+\alpha\frac{p-1}{2p}}\Big( \int r_*^{\gamma\frac{2p}{p-1}}w_R^{-\alpha}\Big)^\frac{p-1}{2p}\nonumber\\
\hspace{-0.7cm}&&\hspace{1cm}\times\Big( \int \calC^4\mu_*^4[\nabla g]_\infty^4\Big)^\frac12\Big( \int w_R^{2p(d-\gamma)+\alpha(p-1)}|f|^{4p}\Big)^\frac1{2p},\label{est2-G2}
\\
\|G_3\|_{\ell}^2&\lesssim_{\gamma,\alpha,p}&\ell^{2d\frac{p-1}{p}}r_*(0)^{2d\frac{p-1}{p}}\Big( \int \calC^p\mu_*^p[\nabla(g\nabla\bar u)]_\infty^p\Big)^\frac{2}p\nonumber\\
\hspace{-0.7cm}&&+\ell^{\gamma}R^{d-\gamma}r_*(0)^{d-\gamma+\alpha\frac{p-1}p}\Big( \int r_*^{\gamma\frac{p}{p-1}}w_R^{-\alpha}\Big)^\frac{p-1}p\nonumber\\
\hspace{-0.7cm}&&\hspace{1cm}\times\Big( \int w_R^{p(d-\gamma)+ \alpha(p-1)}\calC^{2p}\mu_*^{2p}[\nabla(g\nabla\bar u)]_\infty^{2p}\Big)^\frac1p.\label{est2-G3}
\end{eqnarray}
\endgroup
Replacing $p$ by $\frac{2p}{p+1}$, estimating $[\nabla(g\nabla\bar u)]_\infty\lesssim [\nabla g]_\infty[\nabla\bar u]_\infty+[g]_\infty[\nabla^2\bar u]_\infty$, and using H\"older's inequality, the estimate~\eqref{est2-G3} easily leads to
\begin{align*}
\|G_3\|_{\ell}^2&\,\lesssim_{\gamma,\alpha,p}\,\ell^{d\frac{p-1}{p}}r_*(0)^{d\frac{p-1}{p}}\bigg(\Big( \int \calC^2\mu_*^2[\nabla g]_\infty^{2}\Big)\Big( \int [\nabla\bar u]_\infty^{2p}\Big)^\frac{1}p\\
&\hspace{5cm}+\Big( \int \calC^2\mu_*^2[\nabla^2\bar u]_\infty^{2}\Big)\Big( \int [g]_\infty^{2p}\Big)^\frac{1}p\bigg)\nonumber\\&+\ell^{\gamma}R^{d-\gamma}r_*(0)^{d-\gamma+\alpha\frac{p-1}{2p}}\Big( \int r_*^{\gamma\frac{2p}{p-1}}w_R^{-\alpha}\Big)^\frac{p-1}{2p}\\
&\qquad\times\bigg(\Big( \int \calC^4\mu_*^4[\nabla g]_\infty^4\Big)^\frac{1}{2}\Big( \int w_R^{2p(d-\gamma)+\alpha(p-1)}[\nabla\bar u]_\infty^{4p}\Big)^\frac{1}{2p}\\
&\qquad\qquad\qquad+\Big( \int \calC^4\mu_*^4[\nabla^2\bar u]_\infty^4\Big)^\frac{1}{2}\Big( \int w_R^{2p(d-\gamma)+\alpha(p-1)}[g]_\infty^{4p}\Big)^\frac{1}{2p}\bigg).
\end{align*}
Starting from~\eqref{e.operator-norm}
and appealing to~\eqref{est1-G1}--\eqref{est1-G3} for $\ell \le R$ and to~\eqref{est2-G1}--\eqref{est2-G3} for $\ell> R$, 
we obtain the desired estimate~\eqref{est2-WLSI} after arguing as in Substep~2.3 of the proof of Proposition~\ref{prop:boundG2}.
The rest of this step is split into three parts and is dedicated to the proof of~\eqref{est2-G1}--\eqref{est2-G3}.

\medskip

\substep{2.1} Proof of~\eqref{est2-G1}.\\
By \eqref{e.int-r_*} in Lemma~\ref{lem:integrals} and Cauchy-Schwarz' inequality,
\begin{multline*} 
\|G_1\|_{\ell}^2\,\lesssim\,\bigg(\int_{|x|\ge \ell} \big(\ell+r_*(x)\big)^{2d}\Big(\fint_{B_{*}(x)}|g|^2|\nabla\phi+\Id|^2\Big)^2 dx\bigg)^\frac12\\
\times\bigg( \int\Big(\fint_{B_{2\ell*}(x)} (|\nabla w_f|+|\phi||\nabla^2\bar u|)^2\Big)^2dx\bigg)^\frac12
\\
+\bigg(\int_{B_{7\ell*}(0)} \Big(\fint_{\bar B_*(x)}|g|^2|\nabla\phi+\Id|^2\Big)\,dx\bigg)\bigg( \int_{B_{7\ell*}(0)} \Big(\fint_{\bar B_*(x)}(|\nabla w_f|+|\phi||\nabla^2\bar u|)^2\Big)\,dx\bigg).
\end{multline*}
First, we take local suprema of $g$, apply Lemma~\ref{lem:bd-sigma} to control correctors, and use~\eqref{e.equiv-B*}, to the effect of
\begin{align*}
\int_{B_{7\ell*}(0)}\Big(\fint_{\bar B_*(x)}|g|^2|\nabla\phi+\Id|^2\Big)\,dx\,\lesssim\, \int r_*^d[g]_\infty^2.
\end{align*}
Second, using~\eqref{e.equiv-B*} and the energy estimate for~\eqref{eq:wf}, taking local suprema of $\nabla^2\bar u$, and using Lemma~\ref{lem:bd-sigma} to control correctors, we find
\begin{multline*}
\int_{B_{7\ell*}(0)}\Big(\fint_{\bar B_*(x)}(|\nabla w_f|+|\phi||\nabla^2\bar u|)^2\Big)\,dx \,\lesssim \; \int (|\nabla w_f|+|\phi||\nabla^2\bar u|)^2\\
\,\lesssim\, \int (|\phi|+|\sigma|)^2|\nabla^2\bar u|^2\,\lesssim\, \int \calC^2\mu_*^2[\nabla^2\bar u]_\infty^2.
\end{multline*}
Third, appealing to~\eqref{eq-00-pre+} with $s=q=2$ and $|h| =|g||\nabla\phi+\Id|$, we obtain for all $0\le \gamma\le d$, $\alpha>d$, and $p>1$ with $2p(d-\gamma)+\alpha(p-1)<d(4p-1)$,
\begin{multline*}
\int_{|x|\ge \ell} \big(\ell+r_*(x)\big)^{2d}\Big(\int_{B_{*}(x)}|g|^2|\nabla\phi+\Id|^2\Big)^2dx \,
\lesssim_{\gamma,\alpha,p}\,\ell^{2\gamma}R^{2(d-\gamma)}r_*(0)^{2(d-\gamma)+\alpha\frac{p-1}p}\\
\times\Big( \int r_*^{\gamma\frac{2p}{p-1}}w_R^{-\alpha}\Big)^\frac{p-1}p\Big( \int w_R^{2p(d-\gamma)+ \alpha(p-1)}[g(\nabla\phi+\Id)]_2^{4p}\Big)^\frac1p,
\end{multline*}
while the mean-value property~\eqref{eq:mean-value-corr} yields $[g(\nabla\phi+\Id)]_2^{4p}\lesssim r_*^{2pd}[g]_\infty^{4p}$.
The conclusion~\eqref{est2-G1} then follows from the combination of the above estimates with~\eqref{est1-G3-02} (with~$r_*$ replaced by~$2\ell+r_*$). 

\medskip
\substep{2.2} Proof of~\eqref{est2-G2}.

\noindent By \eqref{e.int-r_*} in Lemma~\ref{lem:integrals} and Cauchy-Schwarz' inequality,
\begin{multline*}
\|G_2\|_{\ell}^2\,\lesssim\, \bigg(\int_{|x|\ge \ell}\big(\ell+r_*(x)\big)^{2d}\Big(\fint_{B_{*}(x)}|\nabla u|^2\Big)^2 dx\bigg)^\frac12\\
\times\bigg( \int\Big(\fint_{B_{2\ell*}(x)} (|\phi||\nabla g|+|\nabla r|)^2\Big)^2dx\bigg)^\frac12
\\
+\bigg(\int_{B_{7\ell*}(0)} \Big(\fint_{\bar B_*(x)}|\nabla u|^2\Big)\,dx\bigg)\bigg( \int_{B_{7\ell*}(0)}\Big(\fint_{\bar B_*(x)}(|\phi||\nabla g|+|\nabla r|)^2\Big)\,dx\bigg).
\end{multline*}
First, using~\eqref{e.equiv-B*}, the energy estimate for~\eqref{eq:aux-reps}, taking local suprema of $g$, and applying Lemma~\ref{lem:bd-sigma} to control correctors,
\begingroup\allowdisplaybreaks
\begin{multline*}
\int_{B_{7\ell*}(0)}\Big( \fint_{\bar B_*(x)}(|\phi||\nabla g|+|\nabla r|)^2\Big)\,dx \,\lesssim\, \int (|\phi||\nabla g|+|\nabla r|)^2
\\
\,\lesssim\,\int (|\phi|+|\sigma|)^2|\nabla g|^2\,\lesssim\,\int \calC^2\mu_*^2[\nabla g]_\infty^2.
\end{multline*}
\endgroup
Second, the energy estimate for~\eqref{eq:first-def-ups} yields
\begin{eqnarray*}
\int_{B_{7\ell*}(0)}\Big( \fint_{\bar B_*(x)}|\nabla u|^2\Big)\,dx\,\lesssim\, \int|\nabla u|^2 \,\lesssim\, \int |f|^{2}.
\end{eqnarray*}
Third, appealing to~\eqref{eq-00-pre+} with $s=q=2$ and $h=u$, we obtain for all $0\le\gamma\le d$, $\alpha>d$, and $p>1$ with $2p(d-\gamma)+\alpha(p-1)<d(4p-1)$,
\begin{multline*} 
\int_{|x|\ge \ell} \big(\ell+r_*(x)\big)^{2d} \Big(\fint_{B_{*}(x)}|\nabla u|^2\Big)^2 dx\,\lesssim_{\gamma,\alpha,p}\,\ell^{2\gamma}R^{2(d-\gamma)}r_*(0)^{2(d-\gamma)+\alpha\frac{p-1}p}\\
\times\Big( \int r_*^{\gamma\frac{2p}{p-1}}w_R^{-\alpha}\Big)^\frac{p-1}p
\Big( \int w_R^{2p(d-\gamma)+\alpha(p-1)}|f|^{4p}\Big)^\frac1p.
\end{multline*}
Finally, applying~\eqref{eq-00-pre-1} with $q=2$ (with $r_*$ replaced by $2\ell+r_*$) to the solution $r$ of~\eqref{eq:aux-reps}, taking local suprema of $g$, and applying Lemma~\ref{lem:bd-sigma} to control correctors,
\begin{eqnarray*}
\int \Big(\fint_{B_{2\ell*}(x)}(|\phi||\nabla g|+|\nabla r|)^2\Big)^2dx
\,\lesssim\, \int [(|\phi|+|\sigma|)\nabla g]_2^4\,\lesssim\, \int \calC^4\mu_*^4[\nabla g]_\infty^4.
\end{eqnarray*}
The combination of these four estimates yields the conclusion~\eqref{est2-G2}.

\medskip

\substep{2.3} Proof of~\eqref{est2-G3}.

\nopagebreak

\noindent By \eqref{e.int-r_*} in Lemma~\ref{lem:integrals}, Cauchy-Schwarz' inequality, and the mean-value property~\eqref{eq:mean-value-corr},
\begin{multline}\label{est2-G3-00}
\|G_3\|_{\ell}^2\,\lesssim\,\int_{|x|\ge \ell}\big( \ell+r_*(x)\big)^d \Big(\fint_{B_*(x)} \big(|\phi^*||\nabla(g\nabla\bar u)|+|\nabla R|\big)^2\Big)\,dx\\
+\bigg(\int_{B_{7\ell*}(0)} \Big(\fint_{\bar B_*(x)}\big(|\phi^*||\nabla(g\nabla\bar u)|+|\nabla R|\big)^2\Big)^\frac12dx\bigg)^2,
\end{multline}
We start with the second RHS term. By H\"older's inequality with exponents~$(\frac p{p-1},p)$,
\begin{eqnarray*}
\lefteqn{\int_{B_{7\ell*}(0)} \Big(\fint_{\bar B_*(x)}\big(|\phi^*||\nabla(g\nabla\bar u)|+|\nabla R|\big)^2\Big)^\frac12dx}\\
&\lesssim & \ell^{d\frac{p-1}p}r_*(0)^{d\frac{p-1}p}\bigg(\int_{B_{7\ell*}(0)} \Big(\fint_{\bar B_*(x)}\big(|\phi^*||\nabla(g\nabla\bar u)|+|\nabla R|\big)^2\Big)^\frac p2dx\bigg)^\frac 1p \\
&\lesssim& \ell^{d\frac{p-1}p}r_*(0)^{d\frac{p-1}p}\bigg(\int \Big(\fint_{B_*(x)}\big(|\phi^*||\nabla(g\nabla\bar u)|+|\nabla R|\big)^2\Big)^\frac p2dx\bigg)^\frac 1p,
\end{eqnarray*}
where we passed from $\bar B_*(x)=B_{5r_*(x)}(x)$ to $B_*(x)=B_{r_*(x)}(x)$.
Appealing to the large-scale Cader\'on-Zygmund estimate of Lemma~\ref{lem:r*}(b) with exponent $1<p\le2$ applied to the solution $R$ of~\eqref{eq:aux-Reps}, taking local suprema of $\nabla(g\nabla\bar u)$, and using Lemma~\ref{lem:bd-sigma} to control correctors, we deduce
\begin{equation*}
 \int_{B_{7\ell*}(0)} \Big(\fint_{\bar B_*(x)}\big(|\phi^*||\nabla(g\nabla\bar u)|+|\nabla R|\big)^2\Big)^\frac12dx
\,\lesssim\, \ell^{d\frac{p-1}p}r_*(0)^{d\frac{p-1}p}\Big(\int \calC^p\mu_*^p[\nabla(g\nabla\bar u)]_\infty^p\Big)^\frac 1p.
\end{equation*}
We turn to the first RHS term in~\eqref{est2-G3-00}, apply~\eqref{eq-00-pre+} with $s=q=1$ to~\eqref{eq:aux-Reps}, take local suprema of $\nabla(g\nabla\bar u)$, and use Lemma~\ref{lem:bd-sigma} to control correctors. 
The conclusion~\eqref{est2-G3} follows.


\section{Proof of the main results}

We mainly focus on the proof of the statements for the standard LSI ($\beta>d$), and quickly argue how to adapt the argument to general multiscale LSI ($\beta\le d$) in the last step.

\medskip

\step{1} Proof of Proposition~\ref{prop:scaling} for standard LSI ($\beta>d$).

\noindent Let $\F\in C^\infty_c(\R^d)^{d\times d}$.
Starting point is \eqref{est1-LSI} in Proposition~\ref{prop:boundG2}.
By H\"older's inequality, the triangle inequality in probability, and the stationarity of $r_*$, we obtain for all $R\ge1$, $0<\alpha-d\ll1$, $0<p-1\ll_\alpha1$, and $q\gg\frac1{p-1}$,
\begin{align*}
\expec{\Norm{\partial^\fun J_0(\F)}^{2q}}^\frac1q\,\lesssim_{\alpha,p}\,
\expec{\Big(r_*^{d+\alpha\frac{p-1}p}\Big)^q}^\frac1qR^{d\frac{p-1}p}\Big( \int w_R^{\alpha(p-1)}|F|^{2p}\Big)^\frac1p.
\end{align*}
Replacing $\F$ by $\e^\frac{d}{2}\F(\e\cdot)$ and choosing $R=\frac1\e$, this yields
\begin{align}\label{e.th-1.1}
\expec{\Norm{\partial^\fun\widehat J_0^\e(\F)}^{2q}}^\frac1q\,\lesssim_{\alpha,p}\,
\expec{\Big(r_*^{d+\alpha\frac{p-1}p}\Big)^q}^\frac1q\Big( \int w_1^{\alpha(p-1)}|F|^{2p}\Big)^\frac1p.
\end{align}
We now recall the following implication (which follows from multiscale LSI in form of the moment bounds in~\cite[Proposition~3.1(i)]{DG1}; see also \cite[Step~1 of the proof of Theorem~1]{GNO-quant}): for all random variables $Y_1,Y_2$, given $q_0\ge1$ and $\kappa>0$,
\begin{multline}\label{LSI-int}
\expec{\Norm{\partial^\fun Y_1}^{2q}_\beta}^\frac1q \,\le\, \expec{Y_2^q}^\frac1q~~\text{for all $q\ge q_0$}, \quad\text{and}\quad \expec{\exp(Y_2^{\kappa})}\le2\\
\implies\quad \exists\, C\simeq_{q_0,\kappa}1:~ \expec{\exp\Big(\frac1CY_1^{2\frac{\kappa}{1+\kappa}}\Big)}\le 2.
\end{multline}
Using this property and the moment bound of Lemma~\ref{lem:r*} for $r_*$, the estimate~\eqref{e.th-1.1} leads to the conclusion~\eqref{t1.var}.

\medskip

\step{2} Proof of Theorem~\ref{th:main} for standard LSI ($\beta>d$).
\nopagebreak

\noindent Let $f,g\in C^\infty_c(\R^d)$. We split the proof into two substeps: we first improve~\eqref{est2-LSI} to avoid local suprema in the estimate, and then turn to~\eqref{t1.path} itself.

\medskip
\substep{2.1}
Improvement of~\eqref{est2-LSI}: for all~$R\ge1$, $0<\alpha-d\ll1$, and $0<p-1\ll_\alpha1$,
\begin{multline}\label{eq:bounder-0}
\Norm{\partial^\fun E(f,g)}^2\,\lesssim_{\alpha,p}\,
r_*(0)^{\alpha\frac{p-1}{2p}}\Big( \int r_*^{2d\frac{2p}{p-1}}w_R^{-\alpha}\Big)^\frac{p-1}{2p}\\
\times\bigg(\Big( \int \calC^{4}\mu_*^{4}(|\nabla f|+|\nabla^2\bar u|)^{4}\Big)^\frac{1}{2}\Big( \int w_R^{\alpha(p-1)}|g|^{4p}\Big)^\frac1{2p}\\
+\Big( \int \calC^4\mu_*^4|\nabla g|^4\Big)^\frac12\Big( \int w_R^{\alpha(p-1)}(|f|+|\nabla\bar u|)^{4p}\Big)^\frac1{2p}\bigg).
\end{multline}
We first apply~\eqref{est2-LSI} to the averaged functions $f_1$ and $g_1$ defined by $f_1(x):=\fint_{B(x)}f$ and $g_1(x):=\fint_{B(x)}g$. Noting that $[f_1]_\infty\lesssim\fint_{B_2(x)}|f|$ and that the solution $\bar u_1$ of the homogenized equation~\eqref{eq:first-def-ubar-intro} with averaged RHS $f_1$ is given by $\bar u_1(x)=\fint_{B(x)}\bar u$, and using the Lipschitz continuity of $\calC$, we obtain for all $0<\alpha-d\ll1$ and $0<p-1\ll_\alpha1$,
\begin{multline}\label{eq:bounder-0appl}
\Norm{\partial^\fun E(f,g)}^2\,\lesssim_{\alpha,p}\,\Norm{\partial^\fun E(f_1-f,g_1)}^2+\Norm{\partial^\fun E(f,g_1-g)}^2\\
+r_*(0)^{\alpha\frac{p-1}{2p}}\Big( \int r_*^{2d\frac{2p}{p-1}}w_R^{-\alpha}\Big)^\frac{p-1}{2p}\bigg(\Big( \int \calC^{4}\mu_*^{4}|\nabla^2\bar u|^{4}\Big)^\frac{1}{2}\Big( \int w_R^{\alpha(p-1)}|g|^{4p}\Big)^\frac1{2p}\\
+\Big( \int \calC^4\mu_*^4|\nabla g|^4\Big)^\frac12\Big( \int w_R^{\alpha(p-1)}(|f|+|\nabla\bar u|)^{4p}\Big)^\frac1{2p}\bigg).
\end{multline}
It remains to estimate the first two RHS terms of \eqref{eq:bounder-0appl}, which 
we will prove to be small not because the two-scale expansion is accurate, but because $f_1-f$ and $g_1-g$ are small themselves after rescaling.
Arguing as in the proof of~\eqref{eq:decomp-der-E}, we have the alternative formula
\begin{align*}
\partial^\fun E(f,g)=g\otimes\nabla u-g\otimes\nabla\bar u-g\nabla_i\bar u\otimes(\nabla\phi_i+\ee_i)+\nabla t\otimes\nabla u-\nabla T_i\otimes(\nabla\phi_i+\ee_i),
\end{align*}
where the auxiliary fields $t$ and $T$ are the Lax-Milgram solutions in $\R^d$ of
\begin{eqnarray*}
-\nabla\cdot\Aa^*\nabla t&=&\nabla\cdot((\Aa^*-\bar\Aa^*)g),\\
-\nabla\cdot\Aa^*\nabla T_i&=&\nabla\cdot((\Aa^*-\bar\Aa^*)g\nabla_i\bar u).
\end{eqnarray*}
Using this decomposition and arguing as in the proof of Proposition~\ref{prop:boundG2}, we obtain for all $R\ge1$, $0<\alpha-d\ll1$, and $0<p-1\ll_\alpha1$,
\begin{multline*}
\Norm{\partial^\fun E(f,g)}^2\lesssim r_*(0)^{\alpha\frac{p-1}{2p}}\Big( \int r_*^{d\frac{2p}{p-1}}w_R^{-\alpha}\Big)^\frac{p-1}{2p}\\
\times\min\bigg\{\Big( \int (|f|+|\nabla\bar u|)^4\Big)^\frac12\Big( \int w_R^{\alpha(p-1)}|g|^{4p}\Big)^\frac1{2p}\,;\\
\Big( \int |g|^4\Big)^\frac12\Big( \int w_R^{\alpha(p-1)}(|f|+|\nabla\bar u|)^{4p}\Big)^\frac1{2p}\bigg\}.
\end{multline*}
Using in addition that $|f-f_1|\le\int_0^1\fint_{tB}|\nabla f(\cdot+y)|dydt$, this turns into
\begin{multline*}
\Norm{\partial^\fun E(f_1-f,g_1)}+\Norm{\partial^\fun E(f,g_1-g)}\\
\lesssim
r_*(0)^{\alpha\frac{p-1}{2p}}\Big( \int r_*^{d\frac{2p}{p-1}}w_R^{-\alpha}\Big)^\frac{p-1}{2p}\bigg(\Big( \int (|\nabla f|+|\nabla^2\bar u|)^4\Big)^\frac12\Big( \int w_R^{\alpha(p-1)}|g|^{4p}\Big)^\frac1{2p}\\
+\Big( \int |\nabla g|^4\Big)^\frac12\Big( \int w_R^{\alpha(p-1)}(|f|+|\nabla\bar u|)^{4p}\Big)^\frac1{2p}\bigg).
\end{multline*}
Combining this with~\eqref{eq:bounder-0appl} leads to the conclusion~\eqref{eq:bounder-0}.

\medskip

\substep{2.2} Conclusion.

\noindent By H\"older's inequality, the triangle inequality in probability, and the stationarity of $r_*$, the estimate~\eqref{eq:bounder-0} leads to the following: for all $R\ge1$, $0<\alpha-d\ll1$, $0<p-1\ll_\alpha1$, and $q\gg\frac1{p-1}$,
\begin{multline*}
\expec{\Norm{\partial^\fun E(f,g)}^{2q}}^\frac1q\,\lesssim_{\alpha,p}\,\expec{\Big(r_*^{2d+\alpha\frac{p-1}{2p}}\calC^2\Big)^q}^\frac1q\\
\times R^{\frac d2(1-\frac{1}{p})}\bigg(\Big( \int \mu_*^{4}(|\nabla f|+|\nabla^2\bar u|)^{4}\Big)^\frac{1}{2}\Big( \int w_R^{\alpha(p-1)}|g|^{4p}\Big)^\frac1{2p}\\
+\Big( \int \mu_*^4|\nabla g|^4\Big)^\frac12\Big( \int w_R^{\alpha(p-1)}(|f|+|\nabla\bar u|)^{4p}\Big)^\frac1{2p}\bigg).
\end{multline*}
We then apply the standard weighted Calder\'on-Zygmund theory to the constant-coefficient equation~\eqref{eq:first-def-ubar-intro} for $\bar u$, and replace $f$ and $g$ by $\e^\frac d4f(\e\cdot)$ and $\e^\frac d4g(\e\cdot)$. For the choice $R=\frac1\e$, and by the bound $\mu_*(\frac\cdot \e)\lesssim \mu_*(\frac1\e)\mu_*(\cdot)$, this implies
\begin{multline}\label{eq:bounder-11}
\expec{\Norm{\partial^\fun\widehat E^\e(f,g)}^{2q}}^\frac1q\,\lesssim_{\alpha,p}\,\expec{\Big(r_*^{2d+\alpha\frac{p-1}{2p}}\calC^2\Big)^q}^\frac1q\\
\times \e^{2}\mu_*(\tfrac1\e)^2\bigg(\Big( \int \mu_*^{4}|\nabla f|^{4}\Big)^\frac{1}{2}\Big( \int w_1^{\alpha(p-1)}|g|^{4p}\Big)^\frac1{2p}\\
+\Big( \int \mu_*^4|\nabla g|^4\Big)^\frac12\Big( \int w_1^{\alpha(p-1)}|f|^{4p}\Big)^\frac1{2p}\bigg).
\end{multline}
We now recall the following version of H\"older's inequality: for all random variables $Y_1,Y_2$, given $\kappa_1,\kappa_2>0$,
\begin{multline}\label{holder-exp}
\expec{\exp\big(Y_1^{\kappa_1}\big)}\le2\quad\text{and}\quad \expec{\exp\big(Y_2^{\kappa_2}\big)}\le2\\
\implies\quad \exists C\simeq_{\kappa_1,\kappa_2}1:~ \expec{\exp\Big(\frac1C(Y_1Y_2)^{\frac{\kappa_1\kappa_2}{\kappa_1+\kappa_2}}\Big)}<\infty.
\end{multline}
Using this property, the moment bounds of Lemmas~\ref{lem:r*} and~\ref{lem:bd-sigma} for $r_*$ and $\calC$ yield, for all $\eta>0$, $\expecm{\exp(\frac1{C_\eta}(r_*^{2d}\calC^2)^{\frac14-\eta})}\le2$ for some $C_\eta\simeq_\eta1$.
Combining this with~\eqref{eq:bounder-11}, property~\eqref{LSI-int} yields the conclusion~\eqref{t1.path}.

\medskip

\step3 General multiscale LSI ($\beta\le d$).

\noindent We start with Proposition~\ref{prop:scaling}. By H\"older's inequality, the triangle inequality in probability, and the stationarity of $r_*$, the estimate~\eqref{est1-WLSI} in Proposition~\ref{prop:boundG2} leads to the following: for all $R\ge1$, $0<\gamma<\beta$, $0<\alpha-d\ll1$, $0<p-1\ll_{\gamma,\alpha}1$, and $q\gg\frac1{p-1}$,
\begin{multline*}
\expec{\Norm{\partial^\fun J_0(\F)}_\beta^{2q}}^\frac1{q}\,\lesssim_{\gamma,\alpha,p}\,\expec{\Big(r_*^{d+\alpha\frac{p-1}p}\Big)^q}^\frac1{q}R^{2d}\pi_*(R)^{-1}\\
\times\bigg(R^{-\frac{d}p}\Big( \int w_R^{(d-\gamma)p+\alpha(p-1)}|\F|^{2p}\Big)^\frac1p
+R^{-\frac{2d}p}\Big( \int [F]_2^p\Big)^\frac2p\bigg).
\end{multline*}
Replacing $\F$ by $\e^d\pi_*(\tfrac1\e)^\frac12\F(\e\cdot)$ and choosing $R=\frac1\e$, this yields
\begin{multline*}
\expec{\Norm{\partial^\fun\widehat J_0^\e(\F)}_\beta^{2q}}^\frac1{q}\,\lesssim_{\gamma,\alpha,p}\,\expec{\Big(r_*^{d+\alpha\frac{p-1}p}\Big)^q}^\frac1{q}\\
\times\bigg(\Big( \int w_1^{(d-\gamma)p+\alpha(p-1)}|\F|^{2p}\Big)^\frac1p+\Big( \int [F]_2^p\Big)^\frac2p\bigg).
\end{multline*}
The combination of this estimate with property~\eqref{LSI-int} and with the moment bound of Lemma~\ref{lem:r*} for $r_*$ implies the desired estimate~\eqref{t1.var}.

\medskip

\noindent
We finally turn to Theorem~\ref{th:main}.
Arguing as in Substep~2.1 above, we may get rid of local suprema in the estimate~\eqref{est2-WLSI} in Proposition~\ref{prop:boundG1}. Using then H\"older's inequality, the triangle inequality in probability, and the stationarity of $r_*$, we obtain the following: for all $R\ge1$, $0\le\gamma<\beta$, $0<\alpha-d\ll1$, $0<p-1\ll_\alpha1$, and $q\gg\frac1{p-1}$,
\begin{eqnarray*}
\lefteqn{\expec{\Norm{\partial^\fun E(f,g)}^{2q}}^\frac1q\,\lesssim_{\gamma,\alpha,p}\,\expec{\Big(r_*^{2d+\alpha\frac{p-1}{2p}}\calC^2\Big)^q}^\frac1q}\\
&&\times\Bigg(R^{d+\frac d2(1-\frac{1}{p})}\pi_*(R)^{-1}\bigg(\Big( \int \mu_*^4(|\nabla f|+|\nabla^2\bar u|)^4\Big)^\frac{1}{2}\Big( \int w_R^{2p(d-\gamma)+\alpha(p-1)}|g|^{4p}\Big)^\frac{1}{2p}\\
&&\hspace{4cm}+\Big( \int \mu_*^4|\nabla g|^4\Big)^\frac{1}{2}\Big( \int w_R^{2p(d-\gamma)+\alpha(p-1)}(|f|+|\nabla\bar u|)^{4p}\Big)^\frac{1}{2p}\bigg)\\
&&\qquad+R^{-\beta}\Big( \int \mu_*^2(|\nabla f|+|\nabla^2\bar u|)^2\Big)\bigg(\Big( \int |g|^2\Big)+R^{d-\frac{d}{p}}\Big( \int |g|^{2p}\Big)^\frac{1}p\bigg)\\
&&\qquad+R^{-\beta}\Big( \int \mu_*^2|\nabla g|^2\Big)\bigg(\Big( \int (|f|+|\nabla\bar u|)^2\Big)+R^{d-\frac{d}{p}}\Big( \int (|f|+|\nabla\bar u|)^{2p}\Big)^\frac{1}p\bigg)\Bigg).
\end{eqnarray*}
Since in dimension $d\ge2$ the weights $\mu_*^2$ and $\mu_*^4$ always belong to the Muckenhoupt classes $A_2$ and $A_{4}$, respectively,
we may apply the standard weighted Calder\'on-Zygmund theory to the constant-coefficient equation~\eqref{eq:first-def-ubar-intro} for $\bar u$ in order to simplify the above RHS. Replacing then $f$ and $g$ by $\pi_*(\tfrac1\e)^\frac14\e^\frac d2f(\e\cdot)$ and $\pi_*(\tfrac1\e)^\frac14\e^\frac d2g(\e\cdot)$, choosing $R=\frac1\e$, and using the bound $\mu_*(\frac\cdot \e)\lesssim \mu_*(\frac1\e)\mu_*(\cdot)$,
the conclusion~\eqref{t1.path} follows as in Substep~2.2.


\section*{Acknowledgements}
The work of MD is supported by F.R.S.-FNRS through a Research Fellowship and by the CNRS-Momentum program.
MD and AG acknowledge financial support from the European Research Council under
the European Community's Seventh Framework Programme (FP7/2014-2019 Grant Agreement
QUANTHOM 335410).


\bibliographystyle{plain}
\bibliography{biblio}


\end{document}